\documentclass[11pt,reqno]{amsart}

\usepackage[utf8]{inputenc}
\usepackage[margin=1in]{geometry}
\usepackage{amssymb}
\usepackage{mathtools}
\usepackage{eulervm,charter}   %
\usepackage{xcolor}
\definecolor{linkcol}{rgb}{0,0,0.55}
\definecolor{citecol}{rgb}{0,0.45,0}
\definecolor{urlcol}{rgb}{0.55,0,0}
\usepackage[colorlinks=true,linkcolor=linkcol,citecolor=citecol,urlcolor=urlcol]{hyperref}

\newcommand{\R}{\mathbb{R}}
\newcommand{\Z}{\mathbb{Z}}
\newcommand{\norm}[1]{\left\lVert #1 \right\rVert}
\newcommand{\db}{d_{B}}
\newcommand{\di}{d_{\mathrm{I}}}
\newcommand{\D}[1]{\mathcal{D}_{#1}}
\newcommand{\eps}{\varepsilon}

\newcommand{\Mnm}{\mathcal{M}_2^{n,m}}
\newcommand{\Mdnm}{\mathcal{M}_d^{n,m}}
\newcommand{\dM}{d_{\mathcal{M}}}
\newcommand{\asdim}{\operatorname{asdim}}
\newcommand{\Rk}{\operatorname{Rk}}

\newcommand{\one}{\mathbf{1}}

\newtheorem{theorem}{Theorem}[section]
\newtheorem*{theorem*}{Theorem}
\newtheorem{lemma}[theorem]{Lemma}
\newtheorem{proposition}[theorem]{Proposition}
\newtheorem{corollary}[theorem]{Corollary}
\theoremstyle{definition}

\theoremstyle{remark}
\newtheorem{remark}[theorem]{Remark}

\title{Asymptotic Dimension of Persistence Modules of Bounded Presentation}

\author{Sushovan Majhi}
\address{Data Science Program, George Washington University, Washington D.C., USA}
\email{s.majhi@gwu.edu}

\author{Atish Mitra}
\address{Department of Mathematical Sciences, Montana Technological University, Butte, MT, USA}
\email{amitra@mtech.edu}

\author{\v{Z}iga Virk}
\address{Faculty of Computer and Information Science, University of Ljubljana, Slovenia}
\email{zigavirk@gmail.com}


\begin{document}

\begin{abstract}
We compute the asymptotic dimension of spaces of multiparameter persistence modules of bounded presentation under the interleaving distance.
The $d$-parameter modules with at most $n$ generators and $m$ relations have asymptotic dimension, and Assouad--Nagata dimension, exactly $d(n+m)$, the number of coordinates of their generator and relation grades, over every field; for two parameters this is $2(n+m)$.
The space has continuous moduli, already for two parameters and $n,m\ge2$: over an infinite field it contains, at every scale, infinite uniformly separated families of pairwise non-isomorphic modules with one fibered barcode, so it is not doubling even at a fixed fibered barcode; already for $m\ge1$ it is not coarsely equivalent to a proper metric space.
The proof is an induction on the number of distinct grade coordinates: an interleaving small compared with the gaps between these coordinates moves each of them by at most twice its size and preserves the induced grid representation, so the continuous moduli are invisible at large scale.
Each class of modules at finite distance therefore embeds coarsely into Hilbert space; for two parameters, $n,m\ge2$ and an uncountable field, a class containing such a family admits no coarse embedding into a separable Hilbert space.
Under the two-parameter matching distance, on the space of rank invariants, the asymptotic dimension and the Assouad--Nagata dimension are again $2(n+m)$: slope-one slices pin the grid of grades, and over a grid only finitely many rank invariants occur; on a finite net of $S$ slices the asymptotic dimension is at most $2nS$.
\end{abstract}

\maketitle

\section{Introduction}\label{sec:intro}

Asymptotic dimension, introduced by Gromov~\cite{Gromov1993}, counts the independent directions of a metric space at large scale, and finite asymptotic dimension forces a coarse embedding into Hilbert space~\cite{Roe2003}.
For the spaces that persistent homology produces, the asymptotic dimension was determined in one parameter by Mitra and Virk: the diagrams with at most $n$ points, under the bottleneck distance, have asymptotic dimension $2n$~\cite{Mitra2021,MitraVirk2024corrigendum}, and explicit landmark embeddings of them into Hilbert space, with distortion bounds, are constructed in~\cite{Mitra2024,MajhiMitraVirk-coarse}.
Without the bound on the number of points the space of finite diagrams admits no coarse embedding into Hilbert space~\cite{BubenikWagner2020,Mitra2021} (it was earlier shown to fail Yu's property~A under the Wasserstein distances~\cite{BLPY2019}), and Zava~\cite{Zava2025} charts the same divide for subspaces of the Gromov--Hausdorff space.
This paper determines the asymptotic dimension for multiparameter persistence modules of bounded presentation under the interleaving distance, in any number $d\ge2$ of parameters, and under the matching distance for $d=2$.

Two-parameter modules are not classified by any discrete invariant: their representation type is wild~\cite{BauerBotnanOppermannSteen2020}, and neither the rank invariant nor the fibered barcode determines the isomorphism type~\cite{CarlssonZomorodian2009-multidim}.
The interleaving distance, the largest stable distance over prime fields, by Lesnick's universality theorem~\cite{Lesnick2015}, is NP-hard to compute~\cite{BjerkevikBotnanKerber2020}, and it differs from the matching distance, which compares the barcodes of one-parameter slices~\cite{CerriDiFabioFerriEtAl2013,Landi2018}.
A space of modules needs a complexity bound to have a chance of finite asymptotic dimension, and the size of the presentation is the natural one: bounding the number of bars on every slice, the one-parameter condition, already admits every finite metric space isometrically~\cite{MajhiMitraVirk-coarse}.
That companion paper studies coarse embeddings of spaces of persistence diagrams and modules into Hilbert space; the present paper is self-contained, and its theorems do not appear there.
We write $\Mdnm$ for the $d$-parameter modules with a minimal presentation by at most $n$ generators and $m$ relations, and $\Mnm$ when $d=2$.

\subsection*{Results}
The \emph{position count} of $\Mdnm$ is $d(n+m)$, the number of coordinates of its generator and relation grades; for two parameters it is $2(n+m)$.

\begin{theorem*}[Theorem~\ref{thm:position}]
For all $n\ge1$, $m\ge0$ and $d\ge2$, the asymptotic dimension and the Assouad--Nagata dimension of $(\Mdnm,\di)$ are both $d(n+m)$.
In particular $\asdim(\Mnm,\di)=2(n+m)$.
\end{theorem*}

The lower bound comes from isometric $\ell^\infty$-cubes built from hooks (cyclic modules) and free modules placed far apart, in any number of parameters (Proposition~\ref{prop:lower}).
The upper bound is an induction on the number of distinct grade coordinates, glued at each scale by the saturated unions of Bell and Dranishnikov.
At scale $r$, a module with two coordinates within $16r$ of each other is moved within $16r$, with the same presentation matrix, to a module with fewer coordinates.
A module whose coordinates are all farther apart is a point of a grid together with a representation of the grid, and it is rigid: an interleaving that is small compared with the gaps moves every coordinate by at most twice its size (Lemma~\ref{lem:stablegrid}) and preserves the representation (Lemma~\ref{lem:gaps}).
So every isomorphism type of grid representation is a separate copy of a configuration space, and small interleavings never pass between copies.

The continuous moduli are present already for two parameters.
Over an infinite field there are families of pairwise non-isomorphic modules with one fibered barcode at pairwise interleaving distance exactly $1$ already with two generators and two relations (Proposition~\ref{prop:twogen}), and with four generators the cross-ratio of four lines gives another, at pairwise distance between $1$ and $3$ (Proposition~\ref{prop:crossratio}).
Scaling the grades scales the distance, so for $n,m\ge2$ such families occur at every scale and the space is not doubling even at a fixed fibered barcode.
Theorem~\ref{thm:position} says that none of this is seen at large scale.
Over an uncountable field it does leave one trace: for $n,m\ge2$ a coarse embedding into Hilbert space of a class at finite distance that contains such a family must have a nonseparable target (Remark~\ref{rem:separable}).
The matching distance does not see these families at all, and it too has the position count as its asymptotic dimension.

\begin{theorem*}[Theorem~\ref{thm:matching}]
For all $n\ge1$ and $m\ge0$, the space of rank invariants of $\Mnm$ under the matching distance has asymptotic dimension and Assouad--Nagata dimension $2(n+m)$.
\end{theorem*}

\noindent Slope-one slices pin the grid of grades (Lemma~\ref{lem:slice-rigidity}), and over a fixed grid only finitely many rank invariants occur, which the gluing absorbs at a constant cost.
The isomorphism type of the grid representation is not determined by slices.

Further results give explicit control.
For one generator a module is read from two envelopes of cones on a line, and the matching distance is within a factor three of the interleaving distance (Proposition~\ref{prop:dm-onegen}).
For two generators with at most one relation the slices pin every module away from the walls where the generators become comparable (Lemma~\ref{lem:rigidity}, Proposition~\ref{prop:dm-twoone}).
On a finite net of $S$ slices the sliced matching distance gives asymptotic dimension at most $2nS$, exactly $2n$ for one slice when $m\ge n$ (Proposition~\ref{prop:finitenet}).
The lower bound $2(n+m)$ under either distance is Proposition~\ref{prop:lower}, and for free modules the value $2n$ also follows from Bjerkevik's stability theorem and Kasprowski's theorem on quotients by finite groups (Theorem~\ref{thm:free}).

These spaces lie outside the proper setting in which most of coarse geometry is developed.
For $m\ge1$ the modules at finite distance from $0$ admit no coarse embedding into a proper metric space, and the same holds for the diagrams $\D n$ (Proposition~\ref{prop:notproper}): bounded sets contain infinite uniformly separated families at every scale, as in a regular tree of infinite valence.
The spaces of diagrams with at most $n$ points and of modules of bounded presentation are therefore examples, arising from applications, of finite asymptotic dimension without properness or bounded geometry.
Finite asymptotic dimension gives a coarse embedding into Hilbert space with no properness hypothesis~\cite[Ex.~11.5]{Roe2003},~\cite[Thm.~2.12]{Mitra2021}, so the embedding holds here.
The theorems of coarse index theory that assume properness are discussed in Remark~\ref{rem:proper}.

\subsection*{Organization}
Section~\ref{sec:prelim} fixes notation and records the facts on distances, upset modules and asymptotic dimension that we use.
Section~\ref{sec:wild} exhibits the families with continuous moduli, Section~\ref{sec:lower} proves the lower bounds and treats free modules, and Section~\ref{sec:position} proves the position count.
The remaining sections concern the matching distance: Section~\ref{sec:finitenet} treats finite nets of slices, Sections~\ref{sec:onegen} and~\ref{sec:twoone} treat one generator and two generators with at most one relation with explicit control, Section~\ref{sec:matching} proves the position count for all $n$ and $m$, and Section~\ref{sec:disc} concludes.

\section{Preliminaries}\label{sec:prelim}

Throughout, $k$ is a field, $d\ge2$, and $\R^d$ carries the product order, $x\le y$ when $x_i\le y_i$ for every $i$; we write $\one=(1,\dots,1)$ and $x\vee y$ for the componentwise maximum.
The main case is $d=2$: slices and the matching distance are used only there, and statements that mention $\R^2$ or $\Mnm$ are two-parameter statements.

\subsection{Modules and presentations}\label{sec:modules}
A \emph{module} $M$ assigns a vector space $M(x)$ to each $x\in\R^d$ and a linear map $M(x\le y)\colon M(x)\to M(y)$ to each comparable pair, functorially; morphisms are natural transformations, and kernels and cokernels are taken pointwise.
For $\gamma\in\R^d$ the \emph{free module} $k(-\gamma)$ is $k$ on the orthant $\gamma+\R^d_{\ge0}$ and $0$ elsewhere, with identity maps inside the orthant.
For a closed upset $U\subseteq\R^d$ the \emph{indicator module} $k[U]$ has $k$ at the points of $U$, zero outside, and identity maps between points of $U$; for a finite set $A$ we write $U_A=\bigcup_{a\in A}(a+\R^d_{\ge0})$.
The \emph{shift} of $M$ by $\eps\ge0$ is $M[\eps](x)=M(x+\eps\one)$, and the maps $M(x\le x+\eps\one)$ form the \emph{shift morphism} $M\to M[\eps]$.

A \emph{presentation} of $M$ is a morphism $\varphi\colon F_1\to F_0$ between finite direct sums of free modules with cokernel $M$.
Choosing bases, it is given by \emph{generator grades} $g_1,\dots,g_n$, \emph{relation grades} $h_1,\dots,h_m$, and an $n\times m$ matrix over $k$ whose $(i,j)$ entry vanishes unless $g_i\le h_j$; the $j$-th column is the \emph{$j$-th relation}, an element of $F_0(h_j)$.
Every finitely presented module has a \emph{minimal} presentation, unique up to isomorphism and a direct summand of every presentation~\cite[Def.~7.5, Thm.~7.6, Rem.~7.7]{BotnanLesnick2023}.
Its generator grades are the $a$ with $M(a)$ not spanned by the images of the $M(b)$, $b<a$, counted with multiplicity $\beta_0(a)$, and its relation grades are described in the same way through the submodule of relations; these multisets $\beta_0^M,\beta_1^M$ are the \emph{Betti positions} of $M$, and they are isomorphism invariants.
We write $\Mdnm$ for the set of isomorphism classes of modules over $\R^d$ with a minimal presentation by at most $n$ generators and $m$ relations, and $\Mnm$ when $d=2$.
Beyond all its grades a finitely presented module is constant, of dimension its \emph{eventual dimension}.
We write $\D n$ for the space of persistence diagrams with at most $n$ off-diagonal points under the bottleneck distance $\db$.

\subsection{Distances}\label{sec:distances}
An \emph{$\eps$-interleaving} of $M$ and $N$ is a pair of morphisms $f\colon M\to N[\eps]$ and $g\colon N\to M[\eps]$ such that $g[\eps]\circ f$ and $f[\eps]\circ g$ are the shift morphisms by $2\eps$.
The \emph{interleaving distance} $\di(M,N)$ is the infimum of the $\eps$ for which one exists, and $+\infty$ if there is none; it is infinite between modules of different eventual dimension, since beyond all grades the shift maps are isomorphisms, and scaling all grades by $t>0$ scales it by $t$.

We use the following elementary bound repeatedly.

\begin{lemma}[Moving grades]\label{lem:moving}
Let two presentations over $\R^d$ have the same matrix, both valid (an entry in row $i$ and column $j$ is nonzero only if $g_i\le h_j$ and $g'_i\le h'_j$), with generator grades $g_i,g'_i$ and relation grades $h_j,h'_j$, and let $\delta=\max\bigl(\max_i\norm{g_i-g'_i}_\infty,\max_j\norm{h_j-h'_j}_\infty\bigr)$.
Then their cokernels are $\delta$-interleaved.
\end{lemma}

\begin{proof}
The canonical maps $k(-g_i)\to k(-g'_i)[\delta]$ exist because $g'_i\le g_i+\delta\one$, and they carry the $j$-th relation, at $h_j$, to the $j$-th relation of the second presentation pushed from $h'_j$ up to $h_j+\delta\one\ge h'_j$, which vanishes in the second cokernel.
So they induce a morphism of cokernels, and symmetrically; both composites send each generator to itself $2\delta$ higher, so they are the shift maps.
\end{proof}

For the rest of this subsection $d=2$.
A \emph{normalized line} is a map $L(t)=vt+b$ with $b\in\R^2$, $v\in[1,\infty)^2$ and $\min(v_1,v_2)=1$; we write $\mathcal L$ for the set of these, several $b$ parametrizing one line, which is harmless.
The restriction $M\circ L$ is a one-parameter module, a finite direct sum of intervals~\cite[Thm.~3.4]{BotnanLesnick2023}, whose multiset of intervals is the \emph{barcode} $\mathcal B(M\circ L)$; the assignment $L\mapsto\mathcal B(M\circ L)$ is the \emph{fibered barcode} of $M$, equivalent to its \emph{rank invariant} $\Rk M(x,y)=\operatorname{rank}M(x\le y)$, $x\le y$ (Lemma~\ref{lem:zero-distance}), and incomplete as an invariant of the isomorphism type~\cite{CarlssonZomorodian2009-multidim}.
The \emph{matching distance} is
\[
  \dM(M,N)=\sup_{L\in\mathcal L}\db\bigl(\mathcal B(M\circ L),\mathcal B(N\circ L)\bigr),
\]
where $\db$ is the bottleneck distance; with this normalization of the lines it equals the usual weighted supremum over lines with direction $u$, $\max(u_1,u_2)=1$~\cite{CerriDiFabioFerriEtAl2013,Landi2018}, since reparametrizing such a line by $t\mapsto t/\min(u_1,u_2)$ multiplies every barcode endpoint, hence $\db$, by the weight $\min(u_1,u_2)$; it also equals the supremum of~\cite[Def.~6.11]{BotnanLesnick2023} over all directions $v\in[1,\infty)^2$, since replacing $v$ by $cv$, $c\ge1$, divides $\db$ by $c$; and $\dM\le\di$~\cite{Landi2018}.
It is a pseudometric, and we write $(\overline{\mathcal M}{}_2^{\,n,m},\overline d_{\mathcal M})$ for the metric quotient of $(\Mnm,\dM)$ identifying modules at matching distance $0$.
For a \emph{slice net} $\mathcal L_S=\{L_1,\dots,L_S\}\subset\mathcal L$, the \emph{sliced matching distance} is $\dM^S(M,N)=\max_k\db(\mathcal B(M\circ L_k),\mathcal B(N\circ L_k))\le\dM(M,N)$, and $X_S$ is the metric quotient of $(\Mnm,\dM^S)$ identifying modules at distance $0$.

\subsection{Upset modules}\label{sec:upsets}
The interleaving distance between indicator modules is a containment condition~\cite{MajhiMitraVirk-coarse}; we include the short proof.

\begin{lemma}[Upset modules]\label{lem:upset-interleaving}
For nonempty closed upsets $U,V\subseteq\R^d$, $k[U]$ and $k[V]$ are $\eps$-interleaved exactly when $U+\eps\one\subseteq V$ and $V+\eps\one\subseteq U$.
For finite $A=\{u_i\}$ and $B=\{v_j\}$, $U_A+\eps\one\subseteq U_B$ exactly when every $u_i+\eps\one$ dominates some $v_j$.
For $d=2$, moreover, $\dM(k[U],k[V])=\di(k[U],k[V])$, and this equals the supremum over slope-one lines of the absolute difference of the times at which they enter $U$ and $V$.
\end{lemma}

\begin{proof}
A morphism $k[U]\to k[V][\eps]$ is a scalar on each point of $U\cap(V-\eps\one)$, constant along comparable pairs there, hence on all of this set, which is an upset and so closed under joins.
In an $\eps$-interleaving the composite $k[U]\to k[U][2\eps]$ is the identity on $U$, so the first morphism is nonzero on all of $U$, which forces $U\subseteq V-\eps\one$; symmetrically $V\subseteq U-\eps\one$.
Conversely, if both containments hold, the morphisms equal to $1$ on $U$ and on $V$ form an $\eps$-interleaving.
The statement about finite sets holds because $U_A+\eps\one$ is generated by the points $u_i+\eps\one$.
For the last statement, along the slope-one line through $b$ each indicator module restricts to a single ray $[t_U(b),\infty)$, with $t_U(b)=\inf\{t:b+t\one\in U\}$, so the bottleneck distance there is $|t_U(b)-t_V(b)|$.
If this is at most $\eta$ for every $b$, then for $x\in U$ the line through $x$ has $t_U(x)\le0$, hence $t_V(x)\le\eta$ and $x+\eta\one\in V$; so $U+\eta\one\subseteq V$, symmetrically, and $\di\le\eta$.
With $\dM\le\di$ this gives the equalities.
\end{proof}

\subsection{Asymptotic dimension}\label{sec:asdim}
An \emph{extended metric space} is a metric space whose distance may take the value $+\infty$, such as the spaces of modules; its classes at finite mutual distance are metric spaces.
A family of subsets is \emph{$r$-disjoint} if distinct members are at distance at least $r$, and \emph{uniformly bounded} if their diameters are bounded.
An extended metric space $X$ has \emph{asymptotic dimension} $\asdim X\le K$ if for every $r>0$ it is covered by $K+1$ families of uniformly bounded sets, each family $r$-disjoint~\cite{BellDranishnikov2008}; if the diameters can be taken at most $cr$ with $c$ independent of $r>0$, the \emph{Assouad--Nagata dimension} is at most $K$~\cite{LangSchlichenmaier2005}, and it bounds the asymptotic dimension.
The union of $K+1$ such families is a cover in which every ball of radius $r/2$ meets at most $K+1$ members, so these are, up to the constants, the definitions through multiplicity of~\cite{BellDranishnikov2008,LangSchlichenmaier2005}.
We use three standard facts: the product theorem, $\asdim(X\times Y)\le\asdim X+\asdim Y$ for the maximum metric~\cite{BellDranishnikov2008}; Kasprowski's theorem that the quotient of a proper metric space by an isometric action of a finite group has the asymptotic dimension of the space~\cite{Kasprowski2017}; and the following lower bound.

\begin{lemma}[Cubes]\label{lem:cubes}
If an extended metric space contains, for every $R>0$, an isometric copy of the $\ell^\infty$-cube $[0,R]^k\subset\R^k$, then its asymptotic dimension is at least $k$ (cf.~\cite[Lemma~2.10]{Mitra2021}, \cite[Lemma~5.7]{MajhiMitraVirk-coarse}).
\end{lemma}

\begin{proof}
This is the standard proof that $\asdim\R^k\ge k$~\cite{BellDranishnikov2008}, which uses nothing outside one large cube.
Suppose $\asdim X\le k-1$: there are families $\mathcal U_0,\dots,\mathcal U_{k-1}$ of subsets of $X$, of diameter at most some $D$, whose union covers $X$ and in which distinct members of one family are at distance at least $3$.
Take $R>D+2$ and an isometric copy $Q$ of $([0,R]^k,\ell^\infty)$ in $X$.
For $U\in\mathcal U_j$ meeting $Q$ put $U'=\{x\in Q:d(x,U)\le1\}$, a closed subset of $Q$ of diameter at most $D+2<R$, which therefore meets no two opposite faces of the cube.
Distinct members of one family give sets $U'$ at distance at least $1$, hence disjoint, and only finitely many members of each family meet $Q$, since $Q$ is compact and a point chosen in each would form an infinite $1$-separated set.
So the $U'$ form a finite closed cover of $Q$ by sets none of which meets two opposite faces, and by Lebesgue's covering theorem~\cite[Thm.~IV~2]{HurewiczWallman1941} some point of $Q$ lies in at least $k+1$ of them. A point lies in at most one $U'$ per family, hence in at most $k$, a contradiction.
\end{proof}

The upper bounds glue a part near a space of lower complexity to a part that is followed as points, by the saturated unions of Bell and Dranishnikov; we isolate the step.
It is a simple instance of the Hurewicz-type theorems for the asymptotic dimension~\cite{BellDranishnikov2006} and for the Assouad--Nagata dimension~\cite{BDLM2008}.

\begin{lemma}[Shallow and deep parts]\label{lem:gluing}
Let $X$ and $E$ be extended metric spaces and $Y\subseteq X$, with $\asdim Y\le K$ and $\asdim E\le K$.
Suppose that for every $r>0$ there are $T\ge0$, a cover $X=A\cup B$ such that every point of $A$ has a point of $Y$ at distance at most $T$, and a map $\mu\colon B\to E$ such that $d(x,y)<r$ implies $d_E(\mu x,\mu y)<Cr$, with $C$ independent of $r$, and such that preimages of sets of diameter at most $D$ have diameter at most $\phi(D)<\infty$; the data $T$, $A$, $B$, $\mu$ and $\phi$ may depend on $r$.
Then $\asdim X\le K$.
If moreover $Y$ and $E$ have Assouad--Nagata dimension at most $K$, $T\le tr$ and $\phi(D)\le D$, with $t$ independent of $r$, then $X$ has Assouad--Nagata dimension at most $K$.
\end{lemma}

\begin{proof}
Fix $r>0$.
Cover $E$ by $K+1$ families, each $Cr$-disjoint, of mesh at most $D$; their preimages under $\mu$ are $K+1$ families $\mathcal V^0,\dots,\mathcal V^K$ covering $B$, each $r$-disjoint by the hypothesis on $\mu$, of mesh at most $D_B=\phi(D)$.
Choose $\pi\colon A\to Y$ with $d(x,\pi x)\le T$, put $\rho=4D_B+3r+2T$, and cover $Y$ by $K+1$ families, each $\rho$-disjoint, of mesh at most $D'$; their preimages under $\pi$ are families $\mathcal U^0,\dots,\mathcal U^K$ covering $A$, of mesh at most $D'+2T$, each $(4D_B+3r)$-disjoint.
For each $i$, following Bell and Dranishnikov~\cite{BellDranishnikov2008}, let $\mathcal F^i$ consist of the sets $U\cup\bigcup\{V\in\mathcal V^i:d(V,U)<r\}$, $U\in\mathcal U^i$, and of the members of $\mathcal V^i$ at distance at least $r$ from every member of $\mathcal U^i$.
A member of $\mathcal V^i$ lies within $r$ of at most one member of $\mathcal U^i$, as two of those are more than $2r+D_B$ apart; so two saturated sets are at distance at least $(4D_B+3r)-2(r+D_B)\ge r$, and the other pairs are at distance at least $r$ by construction and by the $r$-disjointness of $\mathcal V^i$.
So each $\mathcal F^i$ is $r$-disjoint, of mesh at most $D'+2T+2r+2D_B$, and the $\mathcal F^i$ cover $X$.
Under the additional hypotheses $D_B\le D\le c_ECr$ and $D'\le c_Y\rho$ for constants $c_E,c_Y$, so the mesh is at most a constant multiple of $r$.
\end{proof}

\section{Continuous moduli at a fixed fibered barcode}\label{sec:wild}

The fibered barcode does not determine a module, and what it fails to see includes continuous moduli: families of pairwise non-isomorphic modules, indexed by the field, have identical slices.
This section records these families and bounds them; the asymptotic dimension is computed in Section~\ref{sec:position}.

\begin{proposition}[Two generators suffice]\label{prop:twogen}
Over any field, let $M_\lambda$, $\lambda\in k\setminus\{0,1\}$, be the cokernel of the map $k(-h_1)\oplus k(-h_2)\to k(-g_1)\oplus k(-g_2)$ with matrix $\left(\begin{smallmatrix}1&1\\1&\lambda\end{smallmatrix}\right)$, where $g_1=(0,2)$, $g_2=(2,0)$, $h_1=(2,3)$, $h_2=(3,2)$.
These modules lie in $\mathcal M_2^{2,2}$, are pairwise non-isomorphic, share one fibered barcode, and satisfy $\di(M_\lambda,M_{\lambda'})=1$ for $\lambda\ne\lambda'$.
The family has $|k|-2$ members, none over $\mathbb F_2$; it is infinite over an infinite field.
\end{proposition}

\begin{proof}
The presentation is minimal, as $h_1,h_2$ are incomparable and differ from $g_1,g_2$, and each relation involves both generators; independence of the relations, from $\lambda\ne1$, makes the eventual dimension $0$.
An isomorphism of modules lifts to an isomorphism of minimal presentations~\cite[Thm.~7.6]{BotnanLesnick2023}; as $g_1,g_2$ are incomparable and so are $h_1,h_2$, the graded automorphisms of the free modules are diagonal, and such an isomorphism replaces the matrix $\left(\begin{smallmatrix}\alpha&\gamma\\\beta&\delta\end{smallmatrix}\right)$ by
\[
  \mathrm{diag}(s,t)\begin{pmatrix}\alpha&\gamma\\\beta&\delta\end{pmatrix}\mathrm{diag}(u,v)^{-1};
\]
so $\alpha\delta/\beta\gamma=\lambda$ is an isomorphism invariant.
At a grade $x$ the module is $k^{\{i:g_i\le x\}}$ modulo the relations present, and every relation present at $x$ has both generators present; the vectors $e_1,e_2,e_1+e_2,e_1+\lambda e_2$ are pairwise independent for $\lambda\notin\{0,1\}$, so every rank of a structure map, and hence the fibered barcode, is independent of $\lambda$.
For the lower bound let $\eps<1$ and suppose $(f,g)$ is an $\eps$-interleaving of $M_\lambda$ and $M_{\lambda'}$.
Lift $f$ to a map $\varphi$ of the free modules, $k(-g_1)\oplus k(-g_2)\to(k(-g_1)\oplus k(-g_2))[\eps]$, carrying relations into relations.
As $g_2\not\le g_1+\eps\one$ and $g_1\not\le g_2+\eps\one$, $\varphi$ is diagonal, $e_i\mapsto a_ie_i$.
Since $h_2\not\le h_1+\eps\one$ and $h_1\not\le h_2+\eps\one$, the relations of $M_{\lambda'}$ at $h_1+\eps\one$ are the multiples of $e_1+e_2$, and at $h_2+\eps\one$ those of $e_1+\lambda'e_2$.
So $a_1e_1+a_2e_2$ is a multiple of $e_1+e_2$ and $a_1e_1+\lambda a_2e_2$ a multiple of $e_1+\lambda'e_2$, which forces $a_1=a_2=0$ as $\lambda\ne\lambda'$.
Then $f=0$, so $gf=0$; this contradicts the fact that $gf$ is the shift, which is nonzero at $g_1$ because $g_1+2\eps\one=(2\eps,2+2\eps)$ lies above no other grade.
So $\di\ge1$. For the upper bound, send each generator of $M_\lambda$ at a grade $x$ to the same generator of $M_{\lambda'}$ at $x+\one$, and symmetrically. This is well defined: a relation of $M_\lambda$ is present at $x$ only if $x\ge(2,3)$ or $x\ge(3,2)$, and then $x+\one\ge(3,3)$, where both relations of $M_{\lambda'}$ are present and $M_{\lambda'}$ vanishes. The composites send each generator to itself two steps up, so they are the shift maps, and $\di(M_\lambda,M_{\lambda'})\le1$.
\end{proof}

\noindent Over a fixed configuration of grades there is no such family when $n\le1$ or $m\le1$: the matrix is then a single row or column, on which the diagonal torus $(k^\times)^n\times(k^\times)^m$ of graded automorphisms already has finitely many orbits, at most one per zero pattern, and the full group of graded automorphisms, which contains it, has no more.
So a continuous modulus over fixed grades first appears for $n,m\ge2$, and Proposition~\ref{prop:twogen} is the smallest possible instance; Proposition~\ref{prop:crossratio} below is a larger one.

\begin{corollary}[Interleaving is not doubling at a fixed fibered barcode]\label{cor:nodoubling}
Over an infinite field and for $n,m\ge2$, for every $t>0$ the space $(\Mnm,\di)$ contains an infinite family of modules with one fibered barcode at pairwise interleaving distance exactly $t$; in particular it is not doubling, even at a fixed fibered barcode.
\end{corollary}

\begin{proof}
Scaling all grades by $t$ scales $\di$ by $t$, so the modules of Proposition~\ref{prop:twogen} with scaled grades form such a family.
In a doubling space the ball of radius $2t$ about one of them, which contains them all, is covered by finitely many balls of radius $t/4$, each of diameter less than $t$ and so containing at most one of them, which is impossible for an infinite family.
\end{proof}

For $m\ge1$ bounded sets of these spaces contain infinite uniformly separated families, and this places them outside the proper setting of coarse geometry.
Recall that a metric space is \emph{proper} if its closed balls are compact, and that a \emph{coarse embedding} $f\colon X\to Y$ comes with nondecreasing functions $\rho_-\le\rho_+$ on $[0,\infty)$, $\rho_-(t)\to\infty$, such that $\rho_-(d(x,y))\le d(fx,fy)\le\rho_+(d(x,y))$.

\begin{proposition}[No properness]\label{prop:notproper}
For $n,m\ge1$ and $d\ge2$, the class of the zero module in $(\Mdnm,\di)$, the modules at finite distance from $0$, admits no coarse embedding into a proper metric space; in particular it is not coarsely equivalent to one.
The same holds for $(\D n,\db)$, $n\ge1$.
\end{proposition}

\begin{proof}
Let $X$ be a metric space with a point $o$ and, for every $t>0$, an infinite set $F_t\subseteq X$ of points pairwise at distance $t$ and within $t$ of $o$.
If $f\colon X\to Y$ is a coarse embedding, choose $t$ with $\rho_-(t)>0$: then $f(F_t)$ is an infinite $\rho_-(t)$-separated subset of the closed ball of radius $\rho_+(t)$ about $f(o)$, which is therefore not compact, and $Y$ is not proper.
In $\D1\subseteq\D n$ the diagrams $\{(x,x+2t)\}$, $x\in t\Z$, are pairwise at distance $t$ and within $t$ of the empty diagram.
For modules take the hooks $H_x=H_{x\one,(x+2t)\one}\in\mathcal M_d^{1,1}$, $x\in t\Z$.
The shift $H_x\to H_x[2\eps]$ vanishes exactly when $\eps\ge t$, so the zero maps form an $\eps$-interleaving of $H_x$ with $0$, and with $H_{x'}$, exactly when $\eps\ge t$; and for $\eps<t$ Lemma~\ref{lem:hooks} applies and would require $|x-x'|\le\eps$.
So $\di(H_x,0)=t$ and $\di(H_x,H_{x'})=t$ for $x\ne x'$; for $d=2$ this is the isometric copy of $\D1$ of Corollary~\ref{cor:onerelation}.
\end{proof}

\noindent The failure of properness of $\D1$ was noted in~\cite{MitraVirk2024corrigendum}, where it prevents a direct application of Kasprowski's theorem.
The regular tree of infinite valence behaves in the same way: it has asymptotic dimension $1$, and the points at distance $t$ from a vertex along distinct edges form such families.
What the families of Proposition~\ref{prop:twogen} add, for $n,m\ge2$ over a field with at least four elements, is that the matching distance vanishes on them, so they are invisible to $\dM$ at every scale.

A second family, with four generators, makes the modulus a cross-ratio, the classical invariant of the quadruples of subspaces of Gelfand and Ponomarev~\cite{GelfandPonomarev1970}.

\begin{proposition}[A cross-ratio family]\label{prop:crossratio}
Over an infinite field there is a family $M_\lambda$, $\lambda\in k\setminus\{0,1\}$, of pairwise non-isomorphic modules in $\mathcal M_2^{4,2}$ with one and the same fibered barcode.
For $\lambda\ne\lambda'$ they satisfy $\dM(M_\lambda,M_{\lambda'})=0$ and $1\le\di(M_\lambda,M_{\lambda'})\le3$.
\end{proposition}

\begin{proof}
Take the four pairwise-incomparable grades $a_1=(0,3)$, $a_2=(1,2)$, $a_3=(2,1)$, $a_4=(3,0)$, for which $a_j-a_i$ has a coordinate at least $1$ whenever $i\ne j$, and $c=(3,3)\ge a_i$.
For $\lambda\in k\setminus\{0,1\}$ let $L_1=\langle e_1\rangle$, $L_2=\langle e_2\rangle$, $L_3=\langle e_1+e_2\rangle$, $L_4=\langle e_1+\lambda e_2\rangle$ in $k^2$, and let $M_\lambda(x)=\sum_{a_i\le x}L_i\subseteq k^2$ with inclusions as structure maps.
Any two of the $L_i$ span $k^2$, so $\dim M_\lambda(x)=\min\{\#\{i:a_i\le x\},2\}$.
The module is generated at the $a_i$; its kernel of presentation is one-dimensional where exactly three generators are present, first at $a_1\vee a_2\vee a_3=(2,3)$ and $a_2\vee a_3\vee a_4=(3,2)$, and two-dimensional at $c$; the relations at $(2,3)$ and $(3,2)$ are independent there, one involving the generator at $a_1$ and the other the generator at $a_4$, so they span it, and the minimal presentation has four generators and two relations, and $M_\lambda\in\mathcal M_2^{4,2}$.
Along any normalized line all structure maps are inclusions, so the restriction has injective structure maps and its barcode consists of essential bars only, as many as the eventual dimension $2$, born at the two smallest entry times of the $a_i$, counted with multiplicity when the slice passes through a join and picks up two generators at once: it is independent of $\lambda$, and $\dM(M_\lambda,M_{\lambda'})=0$.
Moreover $M_\lambda\cong M_{\lambda'}$ only when $\lambda=\lambda'$, the cross-ratio of the four lines being an isomorphism invariant.

We claim $\di(M_\lambda,M_{\lambda'})\ge1$ for $\lambda\ne\lambda'$.
Suppose an $\eps$-interleaving $(f,g)$ exists with $\eps<1$.
Beyond $c$ both modules are the constant $k^2$, so the internal $2\eps$-maps there are the identity and $gf=\mathrm{id}$, $fg=\mathrm{id}$ force $f_c$ to be an automorphism of $k^2$.
At each $a_i$ the shift keeps $a_i+\eps\one\not\ge a_j$ for $j\ne i$, since $a_j-a_i$ has a coordinate at least $1>\eps$, so $M_{\lambda'}(a_i+\eps\one)=L_i^{\lambda'}$ is one-dimensional; naturality of $f$ along $a_i\le c$ gives $f_c(L_i^\lambda)\subseteq L_i^{\lambda'}$, hence equality, for all four $i$.
An automorphism carrying the ordered configuration $(L_i^\lambda)$ to $(L_i^{\lambda'})$ is a projective equivalence of the two, forcing $\lambda=\lambda'$.
For the diameter, $a_i+3\one\ge c$ for every $i$, so $M_{\lambda'}[3]$ is the constant $k^2$ on the support of $M_\lambda$, and the inclusions $M_\lambda(x)\subseteq k^2=M_{\lambda'}(x+3\one)$, in both directions, form a $3$-interleaving.
\end{proof}

The fibers of the position map are at least bounded, once the eventual dimension is fixed; the following lemma is used again only in Remark~\ref{rem:separable}.

\begin{lemma}[Boxes]\label{lem:boxes}
Let $M,N$ be finitely presented modules over $\R^2$ with all generator and relation grades in a box $m_0+[0,w]^2$, and put $J=m_0+w\one$.
If $\dim M(J)=\dim N(J)$ then $\di(M,N)\le w$; otherwise $\di(M,N)=\infty$.
\end{lemma}

\begin{proof}
Beyond $J$ the presentations, and hence $M$ and $N$, have isomorphisms as structure maps, so $M(J)\cong M(y)$ and $N(J)\cong N(y)$ canonically for $y\ge J$.
Choose an isomorphism $\varphi\colon M(J)\to N(J)$.
If $M(x)\ne0$ then $x$ dominates a generator, so $x\ge m_0$ and $x+w\one\ge J$; let $f_x\colon M(x)\to N(x+w\one)$ be the structure map of $M$ to $x+w\one$, followed by the canonical identifications with $M(J)$, by $\varphi$, and by the identification of $N(J)$ with $N(x+w\one)$.
Define $g$ likewise with $\varphi^{-1}$.
Both are natural, being built from structure maps and canonical identifications, and $g[w]\circ f$ and $f[w]\circ g$ are the structure maps over $2w$, so $(f,g)$ is a $w$-interleaving.
Conversely, an $\eps$-interleaving restricted beyond $J$ composes isomorphisms to isomorphisms, so $\dim M(J)=\dim N(J)$.
\end{proof}

\noindent Modules with the same generator and relation grades can differ in eventual dimension, as two generators with two relations killing the same vector or different ones show, so the fibers of the position map must be stratified by $\dim M(J)$; on each stratum they have diameter at most the side of the smallest box containing the grades.
The cross-ratio family is the case $w=3$.

\section{Lower bounds and free modules}\label{sec:lower}

Relations raise the lower bounds, already through modules with a single generator.

\begin{lemma}[Hooks]\label{lem:hooks}
For $a\in\R^d$ and an antichain $B=\{b_1,\dots,b_m\}$, $m\ge0$, with $a\le b_j\ne a$ for every $j$, let $H_{a,B}=k[U_a\setminus U_B]$, $U_a=a+\R^d_{\ge0}$, the cokernel of $\bigoplus_jk(-b_j)\to k(-a)$; its minimal presentation has one generator and the $m$ relations, so $H_{a,B}\in\mathcal M_d^{1,m}$, and every nonzero finitely presented cyclic module is of this form.
If $\eps\ge0$ satisfies $a+2\eps\one\notin U_B$ or $a'+2\eps\one\notin U_{B'}$, then $H_{a,B}$ and $H_{a',B'}$ are $\eps$-interleaved if and only if $\norm{a-a'}_\infty\le\eps$, $U_B+\eps\one\subseteq U_{B'}$ and $U_{B'}+\eps\one\subseteq U_B$.
\end{lemma}

\begin{proof}
Write $H=H_{a,B}$, $H'=H_{a',B'}$.
Suppose $(f,g)$ is an $\eps$-interleaving.
Say $a+2\eps\one\notin U_B$; the other case is symmetric.
The composite $H\to H[2\eps]$ is the shift, nonzero at $a$ because $a+2\eps\one\in U_a\setminus U_B$, so both $f$ and $g$ are nonzero.
As $H$ is generated at $a$, $f$ is determined by $f_a(1)\in H'(a+\eps\one)$, which is therefore nonzero, so $a+\eps\one\ge a'$.
For each $j$, $H(b_j)=0$, so naturality along $a\le b_j$ forces the structure map of $H'$ from $a+\eps\one$ to $b_j+\eps\one$ to kill $f_a(1)$; in an interval module this means $b_j+\eps\one\notin U_{a'}\setminus U_{B'}$, and as $b_j+\eps\one\ge a'$, that $b_j+\eps\one\in U_{B'}$.
So $U_B+\eps\one\subseteq U_{B'}$, and the same argument applied to $g\ne0$ gives $a'+\eps\one\ge a$ and $U_{B'}+\eps\one\subseteq U_B$.
Conversely, given the three conditions, let $f_x=1$ when $x\in U_a\setminus U_B$ and $x+\eps\one\in U_{a'}\setminus U_{B'}$, and $f_x=0$ otherwise, and define $g$ symmetrically.
Since $a'\le a+\eps\one\le x+\eps\one$ for $x\in U_a$, the condition on $x$ is only $x+\eps\one\notin U_{B'}$, and as $U_{B'}$ is an upset it holds at $x$ whenever it holds at some $y\ge x$; where $f_x=1$ and $f_y=0$ the target $H'(y+\eps\one)$ is zero, and a point $y\ge x$ leaving $U_a\setminus U_B$ lies in $U_B$, so $y+\eps\one\in U_{B'}$ and the target map vanishes there too.
So $f$ and $g$ are morphisms.
If $x$ and $x+2\eps\one$ lie in $U_a\setminus U_B$ then $x+\eps\one\notin U_{B'}$, since otherwise $x+2\eps\one\in U_{B'}+\eps\one\subseteq U_B$; so $f_x=g_{x+\eps\one}=1$ and the composites are the shifts.
\end{proof}

\noindent With one relation these are hook modules in the sense of Botnan, Oppermann and Oudot~\cite{Botnan2022signed}, and the interleaving distance between hook modules has an explicit formula due to Botnan, Oppermann, Oudot and Scoccola~\cite[Lemma~4.4]{BOOS2024}; we record it in the present notation, with a proof from Lemma~\ref{lem:hooks}.

\begin{corollary}[One relation]\label{cor:onerelation}
For $a\le b$ in $\R^2$, $a\ne b$, put $p(a,b)=\max_i(b_i-a_i)$.
Then $\di(H_{a,b},0)=\tfrac12p(a,b)$ and
\[
  \di(H_{a,b},H_{a',b'})=\min\Bigl\{\max\bigl(\norm{a-a'}_\infty,\norm{b-b'}_\infty\bigr),\ \tfrac12\max\bigl(p(a,b),p(a',b')\bigr)\Bigr\}.
\]
In particular $(b,d)\mapsto H_{(b,b),(d,d)}$, with the diagonal sent to $0$, is an isometric embedding of $(\D1,\db)$ into $(\mathcal M_2^{1,1},\di)$.
\end{corollary}

\begin{proof}
The shift $H_{a,b}\to H_{a,b}[2\eps]$ vanishes exactly when $x+2\eps\one\ge b$ for every $x\ge a$ with $x\not\ge b$; testing $x=(a_1,y)$ and $x=(y,a_2)$ for large $y$ shows this is $2\eps\ge p(a,b)$.
So the zero maps form an $\eps$-interleaving exactly when $2\eps\ge\max(p,p')$, which gives the distance to $0$ and the bound by $\tfrac12\max(p,p')$.
Below that, $a+2\eps\one\notin U_b$ or $a'+2\eps\one\notin U_{b'}$, and Lemma~\ref{lem:hooks} applies; for single relations $U_b+\eps\one\subseteq U_{b'}$ and its converse together say $\norm{b-b'}_\infty\le\eps$.
For the diagram, $p((b,b),(d,d))=d-b$, so the formula is $\min\{\max(|b-b'|,|d-d'|),\tfrac12\max(d-b,d'-b')\}$, the bottleneck distance on $\D1$.
\end{proof}

\begin{proposition}[Lower bound]\label{prop:lower}
For $n\ge1$ and $m\ge0$,
\[
  \asdim(\Mnm,\di)\ \ge\ 2(n+m),\qquad
  \asdim\bigl(\overline{\mathcal M}{}_2^{\,n,m},\overline d_{\mathcal M}\bigr)\ \ge\ 2(n+m).
\]
Both are realized by direct sums of a hook and free modules, whose presentation matrix has one nonzero entry per column, and on which $\dM$ and $\di$ agree at small scales.
In $d\ge2$ parameters, $\asdim(\mathcal M_d^{n,m},\di)\ge d(n+m)$.
\end{proposition}

\begin{proof}
We begin with one generator, $n=1\le m$.
Fix $R>0$ and consider the hooks $H_{a,B}$ with
\[
  a=(s_0,t_0),\qquad b_j=\bigl((2+j)R,\,(2+m-j)R\bigr)+(s_j,t_j)\quad(j=1,\dots,m),
\]
with $(s_j,t_j)\in[0,R/3]^2$ for $j=0,\dots,m$.
Each $B$ is an antichain above $a$, and for two members at $\ell^\infty$ distance $\eps\le R/3$ we have $a+2\eps\one\le(R,R)$, which lies in no $U_B$, so Lemma~\ref{lem:hooks} applies.
In the containment criterion of Lemma~\ref{lem:upset-interleaving} a shifted corner $b_j+\eps\one$ can dominate only the corner of the same index, since a corner of different index exceeds $b_j$ by at least $2R/3$ in some coordinate; so
\[
  \di(H_{a,B},H_{a',B'})=\max\bigl\{\norm{a-a'}_\infty,\ \max_j\norm{b_j-b_j'}_\infty\bigr\},
\]
the $\ell^\infty$ distance in $\R^{2(m+1)}$.
For $\dM$, every slope-one line $L_b$ meets the support in a bar from $t_{U_a}(b)$ to $t_{U_B}(b)$, of length more than $R$, so for distances below $R/3$ the bottleneck distance on $L_b$ is the larger of $|t_{U_a}(b)-t_{U_{a'}}(b)|$ and $|t_{U_B}(b)-t_{U_{B'}}(b)|$, and the argument of Lemma~\ref{lem:upset-interleaving}, applied to $U_a$ and to $U_B$, gives $\dM=\di$ on the family; distinct members are at positive distance, so the family embeds isometrically in the quotient.
Lemma~\ref{lem:cubes} then gives asymptotic dimension at least $2(m+1)$.

For general $n$ and $m$ take $\Lambda\ge10(m+3)R$ and free generators at $g_i=(-i\Lambda,i\Lambda)+(s'_i,t'_i)$, with $(s'_i,t'_i)\in[0,R/3]^2$.
If $m\ge1$, let $F$ be the free module on $g_1,\dots,g_{n-1}$ and consider the direct sums $F\oplus H$, with $H$ a hook of the family above translated by $(\Lambda,-\Lambda)$; they have $n$ generators and $m$ relations.
If $m=0$, let $F$ be the free module on $g_1,\dots,g_n$ and drop $H$ from what follows.
Every morphism between summands of different kind into an $\eps$-shift, $\eps\le R/3$, vanishes: $k(-g_i)\to H'[\eps]$ because $H'(g_i+\eps\one)=0$, $H\to F'[\eps]$ because $F'(a+\eps\one)=0$ at the generator $a$ of $H$, and $k(-g_i)\to k(-g'_j)[\eps]$ for $j\ne i$ because $g'_j\not\le g_i+\eps\one$.
So every $\eps$-interleaving is block diagonal, and $\di(F\oplus H,F'\oplus H')=\max\{\max_i\norm{g_i-g'_i}_\infty,\di(H,H')\}$, the $\ell^\infty$ distance in $\R^{2(n+m)}$.
For $\dM$, restrict to the slope-one line through $c_u=(u,-u)$: the module $F\oplus H$ has the infinite bars $[t_{g_i}(u),\infty)$, with $t_{g_i}(u)=\max(-i\Lambda+s'_i-u,\,i\Lambda+t'_i+u)$, and the finite bar of $H$, of length more than $R$.
For $u\ge0$ the entry times are $i\Lambda+t'_i+u$, and for $u\le-(n+1)\Lambda$ they are $-i\Lambda+s'_i-u$; on these lines consecutive infinite bars are more than $\Lambda/2$ apart.
So if two members of the family are at matching distance $\vartheta\le R/3$, a matching of cost at most $\vartheta$ on such a line pairs the $i$-th infinite bars, as infinite bars are matched to infinite bars and those of different index are more than $\Lambda/2>\vartheta$ apart, which gives $|t'_i-t''_i|\le\vartheta$ and $|s'_i-s''_i|\le\vartheta$, that is, $\norm{g_i-g'_i}_\infty\le\vartheta$; and on every line it pairs the finite bars, of length more than $R$, so $\sup_u\max(|t_a-t_{a'}|,|t_B-t_{B'}|)\le\vartheta$; by Lemma~\ref{lem:upset-interleaving} this says $\norm{a-a'}_\infty\le\vartheta$, $U_B+\vartheta\one\subseteq U_{B'}$ and $U_{B'}+\vartheta\one\subseteq U_B$, and $\di(H,H')\le\vartheta$ by Lemma~\ref{lem:hooks}.
Hence $\dM=\di$ on the family, distinct members are at positive distance, and Lemma~\ref{lem:cubes} applies under both distances.

In $d$ parameters give every grade $d-2$ further coordinates: those of $a$ and of the $g_i$ in $[0,R/3]$, those of the $b_j$ in $2R+[0,R/3]$.
Lemmas~\ref{lem:upset-interleaving} and~\ref{lem:hooks} hold in $\R^d$, and every comparison above is decided in the first two coordinates, so each $B$ is still an antichain above $a$, the point $a+2\eps\one\le R\one$ lies in no $U_B$, a shifted corner dominates only the corner of the same index, and morphisms between summands of different kind still vanish.
So $\di$ is the $\ell^\infty$ distance in $\R^{d(n+m)}$ on the family, and Lemma~\ref{lem:cubes} gives $d(n+m)$.
\end{proof}

\noindent So relation grades are long features as freely placed as generators.

\begin{theorem}[Free modules]\label{thm:free}
For every $n\ge1$, $\asdim(\mathcal M_2^{n,0},\di)=2n$; in $d\ge2$ parameters the value is $dn$.
\end{theorem}

\begin{proof}
A module with no relations is free, and a free module is determined up to isomorphism by the multiset of grades of its generators; so $\mathcal M^{n,0}_2=\bigsqcup_{r\le n}\mathcal C_r$, where $\mathcal C_r=(\R^2)^r/S_r$ is the space of $r$-point multisets, carrying the interleaving distance between the corresponding free modules $F_A=\bigoplus_{a\in A}k(-a)$.
Pieces of different rank are at infinite distance, so the asymptotic dimension is the maximum over the pieces: the structure maps of a free module are injective, so in an $\eps$-interleaving $(f,g)$ the composite $gf$, a shift of $F_A$, is injective and $f$ is injective at every grade; at a grade $x$ beyond all generators of both modules this gives $k^r\hookrightarrow k^{r'}$ and symmetrically, so $r=r'$.
It therefore suffices to show $\asdim(\mathcal C_r,\di)=2r$.
Write $\varrho(A,B)=\min_\sigma\max_i\norm{a_i-b_{\sigma(i)}}_\infty$ for the $\ell^\infty$-matching distance on $\mathcal C_r$, which is the quotient metric of $(\R^{2r},\norm{\cdot}_\infty)$ by the isometric action of $S_r$ permuting the points.
A matching $\sigma$ attaining $\varrho(A,B)$ pairs the summands into block interleavings, the canonical maps $k(-a_i)\rightleftarrows k(-b_{\sigma(i)})$ at scale $\norm{a_i-b_{\sigma(i)}}_\infty$, so $\di(F_A,F_B)\le \varrho(A,B)$; conversely, the bottleneck distance between free modules is $\varrho(A,B)$, since $\di(k(-a),k(-b))=\norm{a-b}_\infty$ and $\di(k(-a),0)=\infty$, and Bjerkevik's stability theorem for free modules~\cite[Thm.~4.12]{Bjerkevik2021}, with constant $d-1$ in $d\ge2$ parameters, gives $\varrho(A,B)\le\di(F_A,F_B)$.
Thus $\di=\varrho$ on $\mathcal C_r$; in $d$ parameters the two are bi-Lipschitz, which is all we use.
Finally $(\R^{2r},\norm{\cdot}_\infty)$ is proper with $\asdim\R^{2r}=2r$, and by Kasprowski's theorem~\cite{Kasprowski2017} the quotient of a proper metric space by a finite group of isometries has the same asymptotic dimension as the space; so $\asdim(\mathcal C_r,\varrho)=2r$.
Kasprowski's theorem applies directly here, since free modules have no diagonal: $\mathcal C_r$ is the full quotient of $\R^{2r}$ by $S_r$.
The same argument with $\R^{dr}$ gives $dr$ in $d\ge2$ parameters.
\end{proof}

\section{The position count}\label{sec:position}

The upper bound comes from the interleavings themselves, with no reference to slices, and it holds in any number of parameters.
For a finitely presented module $M$ over $\R^d$ and $i=1,\dots,d$ write $X_i(M)$ for the set of $i$-th coordinates of the generator and relation grades of a minimal presentation, $\nu(M)=\sum_i|X_i(M)|\le d(n+m)$ for $M\in\Mdnm$, and $\gamma(M)$ for the smallest distance between distinct elements of one of the sets $X_i(M)$, with $\gamma(M)=\infty$ if there are none.
These depend only on the isomorphism type of $M$, since minimal presentations are unique up to isomorphism.

Below the scale of the gaps between grades, nothing moves at all.

\begin{lemma}[Gaps]\label{lem:gaps}
Let $M,N$ be finitely presented modules over $\R^d$ with all generator and relation grades in a grid $\Gamma=X_1\times\dots\times X_d$, $X_i\subset\R$ finite, and let $\gamma$ be the smallest distance between distinct elements of one of the sets $X_i$.
If $M$ and $N$ are $\eps$-interleaved with $2\eps<\gamma$, then $M\cong N$.
\end{lemma}

\begin{proof}
For $x\in\R^d$ let $\lfloor x\rfloor$ be the largest point of $\Gamma$ below $x$, when one exists.
The presentations are constant on the cells $\{x:\lfloor x\rfloor=p\}$, so $M(\lfloor x\rfloor\le x)$ is an isomorphism, $M(x)=0$ when $\lfloor x\rfloor$ does not exist, and $M$ is determined up to isomorphism by its restriction to $\Gamma$; likewise $N$.
For $p\in\Gamma$ the points $p+\eps\one$ and $p+2\eps\one$ lie in the cell of $p$, since $2\eps<\gamma$, so the structure maps from $p$ to them are isomorphisms.
Put $\bar f_p=N(p\le p+\eps\one)^{-1}f_p\colon M(p)\to N(p)$ and $\bar g_p=M(p\le p+\eps\one)^{-1}g_p\colon N(p)\to M(p)$, natural in $p\in\Gamma$ because $f$, $g$ and the structure maps are.
By naturality $\bar g_p\bar f_p=M(p\le p+2\eps\one)^{-1}g_{p+\eps\one}f_p$, which is the identity since $g_{p+\eps\one}f_p$ is the structure map; symmetrically $\bar f_p\bar g_p=\mathrm{id}$.
So $M|_\Gamma\cong N|_\Gamma$, hence $M\cong N$.
\end{proof}

\begin{lemma}[Stable grid]\label{lem:stablegrid}
Let $M$ and $N$ be finitely presented modules over $\R^d$ that are $\eps$-interleaved, with $0<3\eps<\gamma(M)$.
Then for every $i$, every element of $X_i(M)$ lies within $2\eps$ of an element of $X_i(N)$.
\end{lemma}

\begin{proof}
We treat $i=1$; the other coordinates are symmetric.
Write the points of $\R^d$ as $(x,y)$ with $x\in\R$ and $y\in\R^{d-1}$, and $\one'$ for the vector of ones in $\R^{d-1}$.
As in the proof of Lemma~\ref{lem:gaps}, $M$ is constant on the cells of the grid $X_1(M)\times\dots\times X_d(M)$: a structure map $M(p\le q)$ is an isomorphism when $p$ and $q$ have the same largest grid point below them.

Let $x\in X_1(M)$.
We first find a grade $(x,y)$ of the minimal presentation such that $s\colon M(x-3\eps,y)\to M(x,y)$ is not an isomorphism.
If $x$ is the first coordinate of a generator grade $a=(x,y)$, then $s$ is not surjective, because its image lies in the sum of the images of the $M(b)$, $b<a$, which by minimality does not contain the generator at $a$.
Otherwise $x$ is the first coordinate of a relation grade $h=(x,y)$ and of no generator grade.
Then the free module $F$ of the presentation satisfies $F(x-3\eps,y)\cong F(x,y)$, since $3\eps<\gamma(M)$ and no generator grade has first coordinate in $(x-3\eps,x]$; and for the relation submodule $R\subseteq F$, the image of $R(x-3\eps,y)$ in $R(h)$ lies in the sum of the images of the $R(b)$, $b<h$, which by minimality is a proper subspace of $R(h)$.
So the kernel of $s$, which is $R(h)$ modulo the image of $R(x-3\eps,y)$, is nonzero.

Suppose now that $X_1(N)$ has no element within $2\eps$ of $x$, so that no element lies in $(x-2\eps,x+\eps]$.
Put $A=M(x-3\eps,y)$, $B=M(x,y)$, $A'=N(x-2\eps,y+\eps\one')$, $B'=N(x+\eps,y+\eps\one')$, $A_2=M(x-\eps,y+2\eps\one')$ and $B_2=M(x+2\eps,y+2\eps\one')$, and let $(f,g)$ be the interleaving.
The structure maps $A\to A_2$ and $B\to B_2$ are isomorphisms, because $3\eps<\gamma(M)$ and the coordinates of $(x,y)$ lie in the sets $X_i(M)$, so the corresponding points lie in one cell; hence the structure map $s_2\colon A_2\to B_2$ is conjugate to $s$.
The structure map $\phi\colon A'\to B'$ is an isomorphism, because $N$ is constant in the first coordinate on $[x-2\eps,x+\eps]$ for every value of the other coordinates.
Naturality of $f$ gives $f_B\circ s=\phi\circ f_A$, and $f_A$ is injective since $g_{A'}f_A$ is the isomorphism $A\to A_2$; so $s$ is injective.
Naturality of $g$ gives $s_2\circ g_{A'}=g_{B'}\circ\phi$, and $g_{B'}$ is surjective since $g_{B'}f_B$ is the isomorphism $B\to B_2$; so $s_2$ is surjective, and so is $s$.
Thus $s$ is an isomorphism, a contradiction.
\end{proof}

\begin{theorem}[Position count]\label{thm:position}
For all $n\ge1$, $m\ge0$ and $d\ge2$, the asymptotic dimension and the Assouad--Nagata dimension of $(\Mdnm,\di)$ are both $d(n+m)$.
In particular $\asdim(\Mnm,\di)=2(n+m)$.
\end{theorem}

\begin{proof}
The lower bound is Proposition~\ref{prop:lower}, and it holds for the Assouad--Nagata dimension too, as that is at least the asymptotic dimension.
For $K\ge0$ let $Y_K\subseteq\Mdnm$ consist of the modules with $\nu(M)\le K$; then $Y_{d(n+m)}=\Mdnm$, and $Y_0=\{0\}$.
We show that $Y_K$ has Assouad--Nagata dimension at most $K$, and so asymptotic dimension at most $K$, by induction on $K$, with Lemma~\ref{lem:gluing} for $X=Y_K$ and $Y=Y_{K-1}$; distances here may be infinite, which the lemma allows.

Given $r>0$ put $T=16r$, let $B$ consist of the modules with $\nu(M)=K$ and $\gamma(M)>T$, and let $A=Y_K\setminus B$.

\emph{The shallow part.}
Let $M\in A\setminus Y_{K-1}$, with two consecutive elements $x<x'$ of some $X_i(M)$ at distance at most $T$.
Move every $i$-th coordinate equal to $x'$ in a minimal presentation of $M$ to $x$.
The move is order-preserving on coordinates, so the matrix still defines a map of free modules, and by Lemma~\ref{lem:moving} the new cokernel $M^*$ is within $T$ of $M$.
The new presentation has $n$ generators, $m$ relations and $i$-th coordinates in $X_i(M)\setminus\{x'\}$, the other coordinates being unchanged; a minimal presentation of $M^*$ is a summand of it, so $M^*\in Y_{K-1}$.

\emph{The deep part.}
For $M\in B$ let $q_i=|X_i(M)|$, let $\tau(M)$ be the isomorphism type of the restriction of $M$ to the grid $X_1(M)\times\dots\times X_d(M)$, a representation of the poset $\{1,\dots,q_1\}\times\dots\times\{1,\dots,q_d\}$, and let $\mu(M)\in\R^{q_1}\times\dots\times\R^{q_d}$ list the elements of each $X_i(M)$ in increasing order.
Let $E$ be the disjoint union, over all types $\tau$, of copies of $(\R^{q_1}\times\dots\times\R^{q_d},\norm{\cdot}_\infty)$ with $q_1+\dots+q_d=K$, at infinite distance from one another; then the Assouad--Nagata dimension of $E$ is at most $K$, with the constant of $(\R^K,\norm{\cdot}_\infty)$, whose $r$-disjoint families serve every copy at once, distinct copies being at infinite distance.
Send $M$ to $\mu(M)$ in the copy indexed by $\tau(M)$.

If $M,M'\in B$ and $\di(M,M')<r$, choose an $\eps$-interleaving with $0<\eps<r$.
Since $3\eps<T<\gamma(M),\gamma(M')$, Lemma~\ref{lem:stablegrid} applied in both directions puts every element of $X_i(M)$ within $2\eps$ of $X_i(M')$ and conversely; as consecutive elements are more than $4\eps$ apart on both sides, this is an order-preserving bijection $\psi_i\colon X_i(M')\to X_i(M)$ moving each point by at most $2\eps$, for every $i$.
Let $M''$ be the cokernel of a minimal presentation of $M'$ with every coordinate moved by the $\psi_i$.
By Lemma~\ref{lem:moving}, $M'$ and $M''$ are $2\eps$-interleaved, so composing gives a $3\eps$-interleaving of $M$ and $M''$; both have their grades in the grid $X_1(M)\times\dots\times X_d(M)$, whose gaps exceed $T>6\eps$, so $M\cong M''$ by Lemma~\ref{lem:gaps}, which does not need the moved presentation to be minimal.
As $M''$ is $M'$ moved by an order isomorphism of the grids, $\tau(M')=\tau(M'')=\tau(M)$, and $\norm{\mu(M)-\mu(M')}_\infty\le2\eps<2r$.
Conversely, if $\tau(M)=\tau(M')$ and $\norm{\mu(M)-\mu(M')}_\infty\le D$, move a minimal presentation of $M$ by the order isomorphism of the grids: by Lemma~\ref{lem:moving} the result is within $D$ of $M$, and it is isomorphic to $M'$, as both have grades in the grid of $M'$ and the same restriction to it.
So preimages of sets of diameter $D$ have diameter at most $D$, and the second part of Lemma~\ref{lem:gluing}, with $C=2$, $T=16r$, $\phi(D)=D$, the uniform constant of $E$ and the inductive hypothesis for $Y_{K-1}$, starting from the point $Y_0$, gives Assouad--Nagata dimension at most $K$ for $Y_K$.
\end{proof}

\noindent The constants, such as $T=16r$, are not optimized.
The proof uses two devices at each scale: features below the scale are sent to a space of lower complexity, and features above it are followed as points, the two parts joined by a saturated union.
Both come from Mitra and Virk's proof that $\asdim\D n=2n$~\cite{Mitra2021,MitraVirk2024corrigendum}.
The new ingredient is the rigidity of the deep part, which holds across the continuous moduli.

\begin{remark}[The field and the target]\label{rem:separable}
Let $d=2$. Each of the finitely many classes of $(\Mnm,\di)$ at finite distance, one for each eventual dimension by Lemma~\ref{lem:boxes}, has finite asymptotic dimension and so coarsely embeds into a Hilbert space~\cite[Ex.~11.5]{Roe2003},~\cite[Thm.~2.12]{Mitra2021},~\cite{BellDranishnikov2008}; separability of the target depends on the field.
If $k$ is uncountable and $n,m\ge2$, the modules $M_\lambda$ of Proposition~\ref{prop:twogen}, with all grades scaled by $t>0$, form an uncountable subset of one class at pairwise interleaving distance exactly $t$; a coarse embedding into a separable Hilbert space would carry it, for $t$ large, to an uncountable uniformly separated set, which does not exist.
So that class coarsely embeds into a nonseparable Hilbert space and into no separable one.
If $k$ is countable, every class is separable, being a countable union of images of configuration spaces under the maps of Lemma~\ref{lem:moving}, and the target can be taken to be $\ell^2$: precomposing a coarse embedding with a map that sends each point of the class to a point of a countable dense subset within distance $1$ changes it by a bounded amount and makes its image countable.
The asymptotic dimension is $2(n+m)$ over every field.
\end{remark}

\begin{remark}[What the obvious route does not give]\label{rem:positions}
For $d=2$, a natural route to an upper bound on the asymptotic dimension of $(\Mnm,\di)$ is the map to Betti positions, $M\mapsto(\beta_0^M,\beta_1^M)$, a pair of multisets with at most $n+m$ points of $\R^2$ in all, with a Hurewicz-type theorem~\cite{BellDranishnikov2006} over the finitely many comparability types of such configurations.
That map is not coarsely Lipschitz for $\di$: the hooks $H_{a,a+c_0\one}$ and $H_{a',a'+c_0\one}$ are at interleaving distance at most $c_0/2$ for every $a,a'$, since the zero maps form a $c_0/2$-interleaving (the $c_0$-shift of either module is zero), and their positions are $\|a-a'\|_\infty$ apart.
Features of bounded persistence move freely at bounded cost, exactly as points near the diagonal do in $\D n$, where Mitra and Virk~\cite{Mitra2021} absorb them into a single point.
The proof of Theorem~\ref{thm:position} treats short features the same way, by moving nearby grade coordinates together, and the fibers of the position map over fixed positions, where the continuous moduli live, separately, by Lemmas~\ref{lem:stablegrid} and~\ref{lem:gaps}.
Even for one generator the relation grades are not coarse coordinates: for $B_x=\{(1,10),(x,10-\delta)\}$ with $x>1$, Lemma~\ref{lem:hooks} gives $\di(H_{0,B_x},H_{0,\{(1,10)\}})\le\delta$ for every $x$, because the second corner is nearly redundant and slides along the step at bounded cost.
The relation data of a cyclic module are the upset $U_B$ under the upset-interleaving distance, a directional Hausdorff distance on staircases, close to the hyperspaces of finite subsets studied in~\cite{RadulShukel2009,WeighillYamauchiZava2022}.
\end{remark}

\section{Finite nets of slices}\label{sec:finitenet}

Computations see a module through finitely many slices, as does the embedding COMPLEX~\cite{Majhi2026-complex}, built on a sliced matching distance; this section bounds the asymptotic dimension that such a finite net sees.
The zero set of the matching distance is described by the rank invariant.

\begin{lemma}[Zero matching distance]\label{lem:zero-distance}
For finitely presented $M,N$ the following are equivalent: (1) $\dM(M,N)=0$; (2) $\mathcal B(M\circ L)=\mathcal B(N\circ L)$ for every $L\in\mathcal L$; (3) $\Rk M=\Rk N$.
\end{lemma}

\begin{proof}
(1)$\Leftrightarrow$(2): the slice barcodes are finite, so an optimal matching exists, and cost $0$ forces equality, all bars having positive length.
(2)$\Rightarrow$(3): let $a\le b$.
If $a<b$ in both coordinates, the line through $a$ and $b$ is a normalized line $L$, $a=L(s)$, $b=L(t)$, and $\Rk M(a,b)$ is the number of bars of $\mathcal B(M\circ L)$ containing $[s,t]$; if $a=b$ it is the number containing $s$.
If $a\ne b$ and $a_1=b_1$ (the case $a_2=b_2$ is symmetric), then for $\eps>0$ smaller than the distance from $b_1$ to the next larger first coordinate of a grade of $\beta_0^M\cup\beta_1^M$, the map $M(b)\to M(b+(\eps,0))$ is an isomorphism, since no generator or relation enters the strip $(b_1,b_1+\eps]\times(-\infty,b_2]$; so $\Rk M(a,b)=\Rk M(a,b+(\eps,0))$, the previous case applies, and $\eps$ can be taken small enough for $M$ and $N$ at once.
(3)$\Rightarrow$(2): the barcode of $M\circ L$ is determined by its rank function $(s,t)\mapsto\Rk M(L(s),L(t))$ by the structure theorem.
\end{proof}

\noindent So $\overline{\mathcal M}{}_2^{\,n,m}$ is the space of rank invariants of modules with at most $n$ generators and $m$ relations.
Every slice barcode of a finitely presented module is finite, with bars $[b,d)$, $d\in(b,\infty]$, and a module with at most $n$ generators has at most $n$ bars on every slice, whatever its relations: its restriction to a line is generated by the restrictions of the $n$ generators, and each bar is born at a minimal generator of the restriction.
So the slice barcodes lie in the space $\D{\le n}$ of barcodes with at most $n$ bars, essential bars allowed, under $\db$.
A matching of finite cost pairs essential bars with essential bars, at cost $|b-b'|$, and finite bars with finite bars or the diagonal; so $\D{\le n}$ is the union, at infinite mutual distance, of the pieces $E_a\times\D{n-a}$ with the maximum metric, where $E_a$ is the space of $a$-point multisets of $\R$ under the $\ell^\infty$-matching distance, isometric to the sorted cone $\{x_1\le\dots\le x_a\}\subset(\R^a,\norm{\cdot}_\infty)$.

\begin{proposition}[Finite nets]\label{prop:finitenet}
For every $n,m$ and every slice net $\mathcal L_S$, the quotient $X_S=X_S(n,m)$ of $\Mnm$ satisfies $\asdim(X_S,\dM^S)\le2nS$; if $m\ge n$ and $S=1$, equality holds.
\end{proposition}

\begin{proof}
The map $M\mapsto(\mathcal B(M\circ L_k))_{k\le S}$ embeds $X_S$ isometrically into the max-product of $S$ copies of $(\D{\le n},\db)$.
Each piece $E_a\times\D{n-a}$ has asymptotic dimension at most $a+2(n-a)\le2n$, by the product theorem with $\asdim E_a\le a$ and $\asdim\D j=2j$~\cite{Mitra2021,MitraVirk2024corrigendum}; pieces at infinite distance are unioned trivially, and the product theorem over the $S$ factors gives $2nS$.
For the equality, the module $k(-L_1(b))/k(-L_1(d))$, with one generator and one relation, restricts to the bar $[b,d)$ along $L_1$, and direct sums of at most $n$ of them realize every barcode with at most $n$ finite bars when $m\ge n$; so $(X_1,\dM^1)$ contains an isometric copy of $(\D n,\db)$, of asymptotic dimension $2n$.
\end{proof}

\noindent The proof uses only that every slice barcode has at most $n$ bars, so it holds verbatim for the modules with at most $n$ bars on every slice, which under the continuous distances have infinite asymptotic dimension~\cite[Thm.~5.8]{MajhiMitraVirk-coarse}: a finite net cannot tell bounded presentation from bounded bars.
No single net reaches the continuous matching distance either: for each pair of modules it is computed exactly on a finite set of lines, and the set depends on the pair~\cite{KerberLesnickOudot2019}.

\section{The matching distance for one generator}\label{sec:onegen}

For one generator both distances can be read off two functions on a line, and they agree up to a factor three; so the position count holds under the matching distance.
Parametrize the slope-one lines $L_u(t)=c_u+t\one$ by their intersection $c_u=(u,-u)$ with the antidiagonal, $u\in\R$, and write $\kappa_{\sigma,h}(u)=h+|u-\sigma|$ for the \emph{cone} with vertex $(\sigma,h)$.
For a hook $H=H_{a,B}$ in $\R^2$ put
\[
  t_a(u)=\max_i(a_i-c_{u,i})=\kappa_{\sigma_a,h_a}(u),\qquad t_B(u)=\min_j\max_i(b_{j,i}-c_{u,i})=\min_j\kappa_{\sigma_{b_j},h_{b_j}}(u),
\]
with $\sigma_x=(x_1-x_2)/2$ and $h_x=(x_1+x_2)/2$: the times at which the line through $c_u$ enters $U_a$ and $U_B$, so that $t_a\le t_B$.

\begin{lemma}[Entry functions]\label{lem:entry}
$\norm{a-a'}_\infty=\norm{t_a-t_{a'}}_\infty$, the upset-interleaving distance between $U_B$ and $U_{B'}$ is $\norm{t_B-t_{B'}}_\infty$, and with $p(H)=\sup_u(t_B-t_a)(u)$,
\[
  \di(H,H')=\min\Bigl\{\max\bigl(\norm{t_a-t_{a'}}_\infty,\norm{t_B-t_{B'}}_\infty\bigr),\ \tfrac12\max\bigl(p(H),p(H')\bigr)\Bigr\},\qquad\di(H,0)=\tfrac12p(H).
\]
\end{lemma}

\begin{proof}
The first identity is the computation of $\sup_u|t_a-t_{a'}|$, attained as $u\to\pm\infty$.
For the second, the proof of Lemma~\ref{lem:upset-interleaving} shows that the interleaving distance of $k[U]$ and $k[V]$ is the supremum over slope-one lines of the difference of entry times, and moving the base point along a line shifts both entry times equally, so the antidiagonal suffices.
The shift $H\to H[2\eps]$ vanishes exactly when $a+2\eps\one\in U_B$, that is, $2\eps\ge p_0:=\min_j\max_i(b_{j,i}-a_i)$.
Since $a=c_{\sigma_a}+h_a\one$, one has $(t_B-t_a)(\sigma_a)=p_0$; and for every $u$, with $w=c_u-a$, $t_B(u)=\min_j\max_i(b_{j,i}-a_i-w_i)\le p_0+\max_i(-w_i)=p_0+t_a(u)$.
So $p(H)=p_0$.
Hence the zero maps interleave exactly when $2\eps\ge\max(p,p')$; below that one of the hypotheses of Lemma~\ref{lem:hooks} holds, and the lemma with the first two identities gives the rest.
\end{proof}

\begin{proposition}[Matching distance, one generator]\label{prop:dm-onegen}
On $\mathcal M_2^{1,m}$, $\tfrac13\di\le\dM\le\di$, and the lower bound holds with the supremum in $\dM$ taken over the slope-one lines alone.
Hence $\dM$ is an extended metric on $\mathcal M_2^{1,m}$, infinite between modules of different eventual dimension, and $\asdim(\mathcal M_2^{1,m},\dM)=2(m+1)$ for every $m\ge1$.
\end{proposition}

\begin{proof}
By Lemma~\ref{lem:hooks} the nonzero modules of $\mathcal M_2^{1,m}$ are the free modules $k(-a)$ and the hooks $H_{a,B}$ with $1\le|B|\le m$.
The upper bound is $\dM\le\di$ (Section~\ref{sec:distances}).
For the lower bound we use only the slope-one lines $L_u(t)=c_u+t\one$, which lie in $\mathcal L$.
Along $L_u$ a free module $k(-a)$ restricts to the bar $[t_a(u),\infty)$ and a hook $H_{a,B}$ to the bar $[t_a(u),t_B(u))$, empty when $t_a(u)=t_B(u)$.
Free modules are at distance $\norm{a-a'}_\infty=\sup_u|t_a-t_{a'}|$ under both metrics, a free module and a hook or $0$ are at infinite distance under both, and on the line $L_{\sigma_a}$ the hook $H$ has a bar of length $p(H)$, so $\dM(H,0)\ge\tfrac12p(H)=\di(H,0)$ by Lemma~\ref{lem:entry}.

Let $H=H_{a,B}$ and $H'=H_{a',B'}$ be hooks with $s=\dM(H,H')$.
On each line $L_u$ the bottleneck distance is at most $s$, so either the two bars are matched, with endpoints within $s$, or both have length at most $2s$.
If $\max(p(H),p(H'))\le6s$, then $\di(H,H')\le3s$ by Lemma~\ref{lem:entry}.
Otherwise say $p(H)>6s$.
On $L_{\sigma_a}$ the bar of $H$ has length $p(H)>2s$, so it is matched, and the bar of $H'$ there has length more than $4s$; hence $p(H')>2s$, and on $L_{\sigma_{a'}}$ the bars are matched as well.
So $|t_a-t_{a'}|\le s$ at $u=\sigma_a$ and at $u=\sigma_{a'}$.
The difference $t_a-t_{a'}=h_a-h_{a'}+|u-\sigma_a|-|u-\sigma_{a'}|$ is monotone between $\sigma_a$ and $\sigma_{a'}$ and constant outside, so $\norm{t_a-t_{a'}}_\infty\le s$.
At each $u$ either the bars are matched, and $|t_B-t_{B'}|\le s$, or both are of length at most $2s$, and then
\[
  t_B-t_{B'}=(t_B-t_a)+(t_a-t_{a'})-(t_{B'}-t_{a'})\in[-3s,3s].
\]
So $\norm{t_B-t_{B'}}_\infty\le3s$, and $\di(H,H')\le3s$ by Lemma~\ref{lem:entry}.

Since $\dM\ge\tfrac13\di$, distinct modules are at positive matching distance, so the quotient by $\dM=0$ identifies nothing, and the identity is bi-Lipschitz from $(\mathcal M_2^{1,m},\di)$ to $(\mathcal M_2^{1,m},\dM)$.
Asymptotic dimension is a coarse invariant, so Theorem~\ref{thm:position} with $n=1$ gives the value.
\end{proof}

\noindent With two generators and two relations, over a field with at least four elements, the comparison fails: the modules of Proposition~\ref{prop:twogen} are at matching distance $0$ and interleaving distance $1$.
So continuous moduli are one place where the matching and interleaving distances part; incompleteness of the rank invariant is another (Section~\ref{sec:matching}).

\section{The matching distance for two generators}\label{sec:twoone}

With two generators and at most one relation the slope-one slices still determine a module away from the walls where the generators become comparable, and the position count again holds under the matching distance.
Recall that $x\vee y$ is the componentwise maximum, so that $t_{x\vee y}=\max(t_x,t_y)$.

A \emph{cone envelope} is a minimum $f=\min_j\kappa_j$ of finitely many cones; it is piecewise linear with slopes $\pm1$.
Call $\sigma$ a \emph{local minimum} of $f$ if $f$ has slope $-1$ just left of $\sigma$ and $+1$ just right of it; then $\sigma$ is the vertex of a cone $\kappa_j$ with $h_j=f(\sigma)$.
Order the local minima by $(f(\sigma),\sigma)$ lexicographically; the \emph{depth} of a local minimum $\sigma$ other than the least one is the smallest $H-f(\sigma)$ such that the component of $\{f\le H\}$ containing $\sigma$ contains a smaller local minimum, and the least one has depth $\infty$.
This is the elder rule of sublevel-set persistence.

\begin{lemma}[Valleys]\label{lem:valleys}
Let $f$ and $g$ be cone envelopes with the same number of local minima and let $\eps>0$.
If every finite depth of $f$ exceeds $2\eps$ and $\norm{f-g}_\infty\le\eps$, then the local minima of $f$ and $g$ can be matched so that matched vertices $(\sigma,f(\sigma))$, $(\tau,g(\tau))$ satisfy $|\sigma-\tau|\le2\eps$ and $|f(\sigma)-g(\tau)|\le\eps$.
\end{lemma}

\begin{proof}
On a bounded interval $J\ni\sigma$ whose endpoints are maxima of $f|_J$ and which contains no local minimum of $f$ other than $\sigma$, $f$ equals the cone of $\sigma$: $f$ cannot decrease and then increase again inside $J$ without a further local minimum, so it is monotone on each side of $\sigma$, and having slopes $\pm1$ it is $f(\sigma)+|u-\sigma|$.

Let $\sigma$ be a local minimum of $f$, $h=f(\sigma)$, and $C$ the component of $\{f\le h+2\eps\}$ containing $\sigma$.
If $\sigma'\in C$ were another local minimum, $\sigma$ and $\sigma'$ would lie in one component of $\{f\le h+2\eps\}$, so the larger of the two in the order would have depth at most $2\eps$, a contradiction; so $C=[\sigma-2\eps,\sigma+2\eps]$ by the first paragraph.
On $C$, $g(\sigma)\le h+\eps\le g$ on $\partial C$, so $g$ has a local minimum $\tau$ in the interior of $C$, with $|g(\tau)-h|\le\eps$.
The intervals $C$ for distinct minima of $f$ are disjoint, for the same reason, so this injects the minima of $f$ into those of $g$, a bijection as the counts are equal.
\end{proof}

\begin{lemma}[Classification]\label{lem:class21}
A module with a minimal presentation by two generators and one relation is isomorphic either to $k(-g)\oplus H$ with $H$ a hook with one relation, or to
\[
  N(g_1,g_2,h)=\operatorname{coker}\Bigl(k(-h)\xrightarrow{\ \left(\begin{smallmatrix}1\\1\end{smallmatrix}\right)\ }k(-g_1)\oplus k(-g_2)\Bigr)
\]
with $g_1,g_2$ incomparable and $h\ge g_1\vee g_2$.
In the second case we label the generators so that $g_{1,1}<g_{2,1}$ and $g_{1,2}>g_{2,2}$; the triple $\mu(N)=(g_1,g_2,h)\in\R^6$ is then an isomorphism invariant, and we call
\[
  \eta(N)=\min\bigl(g_{2,1}-g_{1,1},\;g_{1,2}-g_{2,2}\bigr)>0
\]
its \emph{incomparability depth}.
\end{lemma}

\begin{proof}
The relation is $\alpha e_1+\beta e_2$ at a grade $h$.
If $\beta=0$ the module is $H_{g_1,h}\oplus k(-g_2)$, and symmetrically if $\alpha=0$.
Otherwise rescaling the generators makes the relation $e_1+e_2$, and $h\ge g_1\vee g_2$.
If $g_1\le g_2$, the automorphism of $k(-g_1)\oplus k(-g_2)$ fixing $e_1$ and sending $e_2$ to $e_2-e_1$, which exists because $g_1\le g_2$, carries $e_1+e_2$ to $e_2$, and the cokernel is $k(-g_1)\oplus H_{g_2,h}$; symmetrically if $g_2\le g_1$.
The grades of the generators and of the relation of a minimal presentation are isomorphism invariants, being the graded Betti numbers, and the labelling is fixed by the coordinates.
\end{proof}

Along the slope-one line $L_u$ the module $k(-g)\oplus H_{a,b}$ restricts to the bars $[t_g,\infty)$ and $[t_a,t_b)$, and $N(g_1,g_2,h)$ to the cokernel of $k[t_h,\infty)\to k[t_{g_1},\infty)\oplus k[t_{g_2},\infty)$, $1\mapsto(1,1)$, whose bars are
\[
  \bigl[\min(t_{g_1},t_{g_2}),\infty\bigr)\quad\text{and}\quad\bigl[t_{g_1\vee g_2},\,t_h\bigr),
\]
the second empty when $t_h=t_{g_1\vee g_2}$: the older generator survives, and the younger one dies when the relation identifies it with the older.
The lower envelope $\min(t_{g_1},t_{g_2})$ is a cone envelope.
Its two cones cross at height $h_{g_1}+\tfrac12(g_{2,1}-g_{1,1})=h_{g_2}+\tfrac12(g_{1,2}-g_{2,2})$, so it has two local minima, at $\sigma_{g_1}<\sigma_{g_2}$, and its finite depth is $\tfrac12\eta(N)$.

\begin{lemma}[Rigidity]\label{lem:rigidity}
Let $N,N'$ be modules of the second kind in Lemma~\ref{lem:class21}, with $\mu(N)=(g_1,g_2,h)$ and $\mu(N')=(g'_1,g'_2,h')$.
Then $\di(N,N')\le\norm{\mu(N)-\mu(N')}_\infty$.
If moreover $\eps>0$, $\dM(N,N')\le\eps$ and $\eta(N),\eta(N')>8\eps$, then $\norm{g_i-g'_i}_\infty\le3\eps$ for $i=1,2$ and $\norm{h-h'}_\infty\le5\eps$.
\end{lemma}

\begin{proof}
Let $\delta=\norm{\mu(N)-\mu(N')}_\infty$.
The canonical maps $k(-g_i)\to k(-g'_i)[\delta]$ exist because $g'_i\le g_i+\delta\one$, and they carry the relation $e_1+e_2$ at $h$ to $e_1+e_2$ at $h+\delta\one\ge h'$, which vanishes in $N'$; so they induce a morphism $N\to N'[\delta]$, and symmetrically.
Both composites send each generator to itself $2\delta$ higher, so they are the shift maps, and $\di(N,N')\le\delta$.

Now suppose $\dM(N,N')\le\eps$.
On every slope-one line each module has exactly one infinite bar, and infinite bars can only be matched to each other, so the lower envelopes $f=\min(t_{g_1},t_{g_2})$ and $f'=\min(t_{g'_1},t_{g'_2})$ satisfy $\norm{f-f'}_\infty\le\eps$.
Both have two local minima and finite depths greater than $4\eps$, so by Lemma~\ref{lem:valleys} their vertices can be matched with $|\sigma-\tau|\le2\eps$ and heights within $\eps$.
The two vertices of $f$ are $\sigma_{g_2}-\sigma_{g_1}\ge\eta(N)>8\eps$ apart, so the matching pairs $g_1$ with $g'_1$ and $g_2$ with $g'_2$; as $x_1=h_x+\sigma_x$ and $x_2=h_x-\sigma_x$, this gives $\norm{g_i-g'_i}_\infty\le3\eps$, and hence $|t_{g_1\vee g_2}-t_{g'_1\vee g'_2}|\le3\eps$ everywhere.
On each line the finite bars are either matched, so that $|t_h-t_{h'}|\le\eps$, or both of length at most $2\eps$, so that
\[
  t_h-t_{h'}=(t_h-t_{g_1\vee g_2})+(t_{g_1\vee g_2}-t_{g'_1\vee g'_2})-(t_{h'}-t_{g'_1\vee g'_2})\in[-5\eps,5\eps].
\]
Hence $\norm{h-h'}_\infty=\sup_u|t_h-t_{h'}|\le5\eps$ by Lemma~\ref{lem:entry}.
\end{proof}

\noindent Away from the comparability walls the configuration is thus read off the slices, and the interleaving and matching distances agree with the configuration distance up to a factor five at small scales.
Near a wall the configuration is not a coarse coordinate: for $g_1=(0,5)$, $g_2=(10,5-\delta)$ and any $h$, moving $g_2$ to $(10,5)$ costs $\delta$ and, by the proof of Lemma~\ref{lem:class21}, gives $k(-g_1)\oplus H_{(10,5),h}$, or $k(-g_1)$ if $h=(10,5)$.
The shallow part of the next proof absorbs exactly these modules.

\begin{proposition}[Matching distance, two generators]\label{prop:dm-twoone}
For $m\in\{0,1\}$, $\asdim\bigl(\overline{\mathcal M}{}_2^{\,2,m},\overline d_{\mathcal M}\bigr)=2(2+m)$.
\end{proposition}

\begin{proof}
The lower bound is Proposition~\ref{prop:lower}.
On every line the number of infinite bars is the eventual dimension, so classes at finite matching distance have a fixed eventual dimension, and we bound each class.
Lemma~\ref{lem:gluing} applies verbatim to pseudometrics, so we work with $\dM$ on modules.

\emph{Free modules of rank two.}
For $F=k(-g_1)\oplus k(-g_2)$ the slope-one line $L_u$ carries the infinite bars $[t_{g_1}(u),\infty)$ and $[t_{g_2}(u),\infty)$, and the bottleneck distance between two such pairs is the larger of the differences of their minima and of their maxima; so $\dM(F,F')\ge\norm{f-f'}_\infty$ for the lower envelopes $f=\min(t_{g_1},t_{g_2})$ and $f'$.
If $g_1\le g_2$ then $t_{g_1}\le t_{g_2}$, so for two comparable pairs $\dM\ge\max_i\sup_u|t_{g_i}-t_{g'_i}|=\max_i\norm{g_i-g'_i}_\infty\ge\di$ by Lemma~\ref{lem:entry} and Theorem~\ref{thm:free}; hence $\dM=\di$ on the comparable pairs, a subspace $Y$ of $(\mathcal C_2,\di)$, of asymptotic dimension at most $4$.
For an incomparable pair, labelled as in Lemma~\ref{lem:class21}, $f$ is the lower envelope above, with two local minima and finite depth $\tfrac12\eta(F)$, where $\eta(F)=\min(g_{2,1}-g_{1,1},\,g_{1,2}-g_{2,2})$; so if $\eps>0$, $\dM(F,F')\le\eps$ and $\eta(F),\eta(F')>8\eps$, the argument for the generators in the proof of Lemma~\ref{lem:rigidity}, applied to $f$ and $f'$, gives $\norm{g_i-g'_i}_\infty\le3\eps$ for $i=1,2$.
Given $r>0$ put $T=16r$, let $B$ consist of the incomparable pairs with $\eta>T$ and $A$ of the rest.
An incomparable pair in $A$ is within $T$ of $Y$ for $\di$, hence for $\dM$: if $\eta(F)=g_{1,2}-g_{2,2}$, the pair $g_1\le g_2^*=(g_{2,1},g_{1,2})$ is comparable and, by Lemma~\ref{lem:moving}, within $\eta(F)\le T$ of $F$; symmetrically if $\eta(F)=g_{2,1}-g_{1,1}$.
On $B$ let $\mu(F)=(g_1,g_2)\in(\R^4,\norm{\cdot}_\infty)$: then $\dM(F,F')<r$ implies $\norm{\mu(F)-\mu(F')}_\infty<3r$, and $\dM\le\di\le\norm{\mu(F)-\mu(F')}_\infty$; Lemma~\ref{lem:gluing} with $K=4$ and $C=3$ gives asymptotic dimension at most $4$.
Rank at most one gives $\{0\}$ and a copy of $(\R^2,\norm{\cdot}_\infty)$, on which $\dM=\di$.

\emph{One relation.}
Eventual dimension $2$ is the free case above.
Eventual dimension $0$ consists of $0$ and the hooks with one relation, of asymptotic dimension $4$ under $\dM$ by Proposition~\ref{prop:dm-onegen}.
Eventual dimension $1$ consists, by Lemma~\ref{lem:class21}, of the decomposable modules $k(-g)\oplus H$, with $H$ a hook with one relation or $0$, which form a set $Y$, and of the modules $N(g_1,g_2,h)$.

By the Krull--Schmidt theorem~\cite[Thm.~4.2]{BotnanLesnick2023}, $k(-g)\oplus H$ determines $g$ and $H$.
Put \[\rho=\max(\norm{g-g'}_\infty,\di(H,H')).\]
Direct sums of interleavings give $\di\le\rho$.
Conversely, on each slope-one line the infinite bars $[t_g,\infty)$ and $[t_{g'},\infty)$ must be matched to each other, so $\dM(k(-g)\oplus H,k(-g')\oplus H')$ is at least $\norm{g-g'}_\infty=\sup_u|t_g-t_{g'}|$ and at least the supremum over slope-one lines of the bottleneck distance of the finite bars, which is at least $\tfrac13\di(H,H')$ by the proof of Proposition~\ref{prop:dm-onegen}.
So $\tfrac13\rho\le\dM\le\di\le\rho$ on $Y$, which under $\dM$ is therefore bi-Lipschitz to the product of $(\R^2,\norm{\cdot}_\infty)$ with the eventual-dimension-$0$ part of $(\mathcal M_2^{1,1},\di)$, and $\asdim(Y,\dM)\le2+4=6$ by the product theorem~\cite{BellDranishnikov2008} and Theorem~\ref{thm:position}.

Given $r>0$ put $T=16r$, let $B$ consist of the modules $N$ with $\eta(N)>T$, and let $A$ consist of $Y$ and the modules $N$ with $\eta(N)\le T$.
Every $N\in A\setminus Y$ is within $T$ of $Y$ for $\di$, hence for $\dM$: if $\eta(N)=g_{1,2}-g_{2,2}$, replace $g_2$ by $g_2^*=(g_{2,1},g_{1,2})$, which satisfies $g_1\le g_2^*\le g_1\vee g_2\le h$; by Lemma~\ref{lem:moving} the cokernel $N^*$ of the same matrix with $g_2^*$ in place of $g_2$ satisfies $\di(N,N^*)\le\eta(N)$, and $N^*\cong k(-g_1)\oplus H_{g_2^*,h}$ by the proof of Lemma~\ref{lem:class21}, with $H_{g_2^*,h}$ read as $0$ when $h=g_2^*$; symmetrically if $\eta(N)=g_{2,1}-g_{1,1}$.
On $B$ take $\mu$ of Lemma~\ref{lem:class21}, with values in $(\R^6,\norm{\cdot}_\infty)$, of asymptotic dimension $6$.
If $N,N'\in B$ and $\dM(N,N')<r$, choose $\eps$ with $\dM(N,N')\le\eps<r$ and $\eps>0$; then $\eta(N),\eta(N')>16r>8\eps$, so Lemma~\ref{lem:rigidity} gives $\norm{\mu(N)-\mu(N')}_\infty\le5\eps<5r$.
And $\dM\le\di\le\norm{\mu(N)-\mu(N')}_\infty$, so preimages of sets of diameter $D$ have diameter at most $D$.
Lemma~\ref{lem:gluing} with $K=6$ and $C=5$ gives asymptotic dimension at most $6$.
\end{proof}

\noindent The next case differs: with two generators and two relations, over a field with at least four elements, the modules of Proposition~\ref{prop:twogen} share every slice, so slices no longer pin the isomorphism type, and the continuous invariant $\lambda$ is seen only by the interleaving distance.
Slices still pin the grid of grades, and that suffices for the asymptotic dimension (Section~\ref{sec:matching}).

\section{The matching distance}\label{sec:matching}

Under the matching distance the position count holds for all $n$ and $m$.
The grid of grade coordinates is a function of the rank invariant and is read off slope-one slices; over a fixed grid only finitely many rank invariants occur; and a gluing that allows boundedly many fibres costs only a constant.
For $M\in\Mnm$ we use $X_i(M)$, $\nu(M)$ and $\gamma(M)$ as in Section~\ref{sec:position}, and write $\Gamma(M)=X_1(M)\times X_2(M)$.
Two facts from that section: (F1) if $p\le q$ have the same largest point of $\Gamma(M)$ below them, $M(p\le q)$ is an isomorphism, and if $p$ has none then $M(p)=0$ (proof of Lemma~\ref{lem:gaps}); (F2) for $x\in X_1(M)$ there is a grade $(x,y_0)$ of a minimal presentation such that $M((x-\delta,y_0)\le(x,y_0))$ is not an isomorphism for every $0<\delta<\gamma(M)$, and symmetrically for $X_2(M)$ (first paragraph of the proof of Lemma~\ref{lem:stablegrid}, with $3\eps$ replaced by $\delta$: a generator coordinate needs no bound on $\delta$, and a relation coordinate only that no generator has first coordinate in $(x-\delta,x]$).

\begin{lemma}[Grid and rank type]\label{lem:grid-rank}
For $M\in\Mnm$, $X_1(M)$ is the set of $x\in\R$ for which there is $y$ with $M((x-\delta,y)\le(x,y))$ not an isomorphism for all small $\delta>0$, and similarly for $X_2(M)$; so $\Gamma(M)$, $\nu(M)$ and $\gamma(M)$ depend only on $\Rk M$.
Moreover $\Rk M$ is determined by the sorted coordinates $\mu(M)\in\R^{q_1}\times\R^{q_2}$, $q_i=|X_i(M)|$, together with the \emph{rank type} $\tau_{\mathrm{Rk}}(M)$, the function $(a,b)\mapsto\Rk M(z_a,z_b)$ on pairs $a\le b$ of points of $\{1,\dots,q_1\}\times\{1,\dots,q_2\}$, where $z_a$ is the grid point of index $a$; and at most $N(n,m)<\infty$ rank types occur in $\Mnm$.
\end{lemma}

\begin{proof}
Every $x\in X_1(M)$ has the property by (F2).
If $x\notin X_1(M)$, then for small $\delta$ the points $(x-\delta,y)$ and $(x,y)$ either have the same largest grid point below them, and the map is an isomorphism by (F1), or neither has one, and it is $0\to0$.
Whether a structure map is an isomorphism is determined by its rank and the dimensions of its source and target, all values of $\Rk M$.
For $p\le q$ with largest grid points $\lfloor p\rfloor\le\lfloor q\rfloor$ below them, $M(\lfloor q\rfloor\le q)\circ M(\lfloor p\rfloor\le\lfloor q\rfloor)=M(p\le q)\circ M(\lfloor p\rfloor\le p)$ with both outer maps isomorphisms by (F1), so $\Rk M(p,q)=\Rk M(\lfloor p\rfloor,\lfloor q\rfloor)$; and $\Rk M(p,q)=0$ if $\lfloor p\rfloor$ does not exist.
A module with $n$ generators has $\dim M(p)\le n$, so $\tau_{\mathrm{Rk}}(M)$ takes values in $\{0,\dots,n\}$ on at most $(q_1q_2)^2\le(n+m)^4$ pairs.
\end{proof}

\noindent By Lemma~\ref{lem:zero-distance}, $X_i$, $\nu$, $\gamma$, $\mu$ and $\tau_{\mathrm{Rk}}$ are therefore well defined on $\overline{\mathcal M}{}_2^{\,n,m}$.

\begin{lemma}[Slices pin the grid]\label{lem:slice-rigidity}
Let $M,N\in\Mnm$, $\eps>0$ and $T>8(n+m)\eps$.
If $\gamma(M),\gamma(N)>T$ and $\dM(M,N)<\eps$, then for $i=1,2$ there is an order-preserving bijection $X_i(M)\to X_i(N)$ moving every point by less than $\eps$; in particular $\norm{\mu(M)-\mu(N)}_\infty<\eps$.
\end{lemma}

\begin{proof}
Let $x\in X_1(M)$ and take $(x,y_0)$ as in (F2).
As $y_0\in X_2(M)$ and $\gamma(M)>T$, no element of $X_2(M)$ lies in $(y_0,y_0+T]$; for $y$ in this interval the crossing map $M((x-\delta,y)\le(x,y))$ is conjugate to the one at $y_0$ by vertical maps, which are isomorphisms by (F1), so it is not an isomorphism.
The set of $y\in[y_0+T/4,y_0+3T/4]$ within $2\eps$ of $X_2(N)$ has measure at most $4(n+m)\eps<T/2$; choose $y^*$ in the interval outside it, and let $L(t)=(x,y^*)+t\one$, a normalized line.
For $0<|t|<T/4$ the point $L(t)$ has first coordinate in $(x-T/4,x+T/4)\setminus\{x\}$, which meets no element of $X_1(M)$, and second coordinate in $(y_0,y_0+T)$, which meets none of $X_2(M)$.
So by (F1) the module $P=M\circ L$ has isomorphisms as structure maps between parameters in $(-T/4,0)$ and between parameters in $[0,T/4)$, and $P(-\delta\le0)$ is the crossing map at height $y^*$ composed with an isomorphism.
A bar of the finite barcode of $P$ has an endpoint at $t$ exactly when $P(t-\delta\le t)$ fails to be an isomorphism for all small $\delta>0$; so some bar has an endpoint at $0$, and it has length at least $T/4$, no other endpoint lying in $(-T/4,T/4)$.
As $\db(\mathcal B(M\circ L),\mathcal B(N\circ L))<\eps$ and this bar costs at least $T/8>\eps$ to match to the diagonal, it is matched to a bar of $N\circ L$, births to births, deaths to deaths and essential bars to essential bars; so the endpoint corresponding to $0$ is a finite $t'$ with $|t'|<\eps$, and $N(L(t'-\delta)\le L(t'))$ is not an isomorphism for small $\delta$.
By (F1) for $N$ a coordinate of $L(t')=(x+t',y^*+t')$ then lies in the grid of $N$, and by the choice of $y^*$ it is the first: $x+t'\in X_1(N)$.
So every element of $X_1(M)$ is within $\eps$ of $X_1(N)$, and conversely; as both sets have gaps larger than $T>2\eps$, nearest points give mutually inverse order-preserving bijections.
The same argument with the coordinates exchanged treats $X_2$.
\end{proof}

\begin{lemma}[Fibres]\label{lem:fibres}
If $M,N\in\Mnm$ have the same rank type and $\norm{\mu(M)-\mu(N)}_\infty\le D$, then $\dM(M,N)\le D$.
\end{lemma}

\begin{proof}
Let $h_i\colon\R\to\R$ be the increasing piecewise-linear homeomorphism sending the $j$-th element of $X_i(N)$ to the $j$-th element of $X_i(M)$, a translation outside them, so that $|h_i(t)-t|\le D$, and put $h=h_1\times h_2$ and $M'=M\circ h$.
As $k(-g)\circ h=k(-h^{-1}(g))$ and cokernels are computed pointwise, moving every grade $g$ of a minimal presentation of $M$ to $h^{-1}(g)$ gives a minimal presentation of $M'$, precomposition with an order automorphism being an autoequivalence; so $M'\in\Mnm$, with grid $\Gamma(N)$ and the rank type of $M$, which is that of $N$.
By Lemma~\ref{lem:grid-rank}, $\Rk M'=\Rk N$, so $\dM(M',N)=0$ by Lemma~\ref{lem:zero-distance}.
Since $x\le h(x+D\one)$ and $h(x)\le x+D\one$, the structure maps $M(x)\to M(h(x+D\one))=M'(x+D\one)$ and $M'(x)=M(h(x))\to M(x+D\one)$ form a $D$-interleaving of $M$ and $M'$.
So $\dM(M,N)\le\dM(M,M')\le\di(M,M')\le D$.
\end{proof}

\noindent Over a grid the fibres of $\mu$ are thus unions of boundedly many sets of small diameter, and the gluing lemma extends to that case.

\begin{lemma}[Gluing with finite fibres]\label{lem:gluing-finite}
Lemma~\ref{lem:gluing} holds with its hypothesis on preimages replaced by the following: the preimage under $\mu$ of every set of diameter at most $D$ is a union of at most $N$ sets of diameter at most $\phi(D)$, with $N$ independent of $r$.
Under its additional hypotheses, the conclusion on the Assouad--Nagata dimension holds as well.
\end{lemma}

\begin{proof}
In the proof of Lemma~\ref{lem:gluing}, fix for each member $U$ of the covering families of $E$ a decomposition of $\mu^{-1}(U)$ into at most $N$ pieces, and replace $\mu^{-1}(U)$ by its $r$-components, the classes of the equivalence relation generated by lying in pieces at distance less than $r$.
A component is a union of at most $N$ pieces connected through such pairs, so its diameter is at most $N\phi(D)+(N-1)r$.
Distinct components of one preimage are at distance at least $r$ by construction, and components of preimages of distinct members of one family are at distance at least $r$ because that family is $Cr$-disjoint.
So the families $\mathcal V^i$ remain $r$-disjoint, now with mesh $D_B=N\phi(D)+(N-1)r$, which is linear in $r$ under the additional hypotheses, and the rest of the proof is unchanged.
\end{proof}

\begin{theorem}[Position count under the matching distance]\label{thm:matching}
For all $n\ge1$ and $m\ge0$, the asymptotic dimension and the Assouad--Nagata dimension of $\bigl(\overline{\mathcal M}{}_2^{\,n,m},\overline d_{\mathcal M}\bigr)$ are both $2(n+m)$, over every field.
\end{theorem}

\begin{proof}
The lower bound is Proposition~\ref{prop:lower}, and it holds for the Assouad--Nagata dimension, which is at least the asymptotic dimension.
For the upper bound let $Y_K$ consist of the classes with $\nu\le K$, so that $Y_0=\{0\}$ and $Y_{2(n+m)}$ is the whole space; we show that $Y_K$ has Assouad--Nagata dimension at most $K$ by induction on $K$, with Lemma~\ref{lem:gluing-finite} for $X=Y_K$ and $Y=Y_{K-1}$.
Given $r>0$ put $T=16(n+m)r$, let $B$ consist of the classes with $\nu=K$ and $\gamma>T$, and let $A=Y_K\setminus B$.

\emph{The shallow part.}
Let $M$, with a minimal presentation, represent a class in $A\setminus Y_{K-1}$; by Lemma~\ref{lem:grid-rank} its coordinates are those of the class, so two consecutive elements $x<x'$ of some $X_i(M)$ are at distance at most $T$.
As in the proof of Theorem~\ref{thm:position}, moving the coordinates equal to $x'$ to $x$ gives $M^*$ with $\nu(M^*)\le K-1$ and $\dM(M,M^*)\le\di(M,M^*)\le T$, by Lemma~\ref{lem:moving}.

\emph{The deep part.}
Let $E$ be the disjoint union, at infinite mutual distance, of the spaces $(\R^{q_1}\times\R^{q_2},\norm{\cdot}_\infty)$ with $q_1+q_2=K$, of Assouad--Nagata dimension $K$ with a uniform constant, and send a class in $B$ to $\mu$ in the copy indexed by $(q_1,q_2)$.
If two classes in $B$, with representatives $M_1$ and $M_2$, are at matching distance less than $r$, choose $\eps<r$ with $\dM(M_1,M_2)<\eps$; as $T>8(n+m)\eps$, Lemma~\ref{lem:slice-rigidity} gives the same $(q_1,q_2)$ and $\norm{\mu(M_1)-\mu(M_2)}_\infty<r$.
The preimage of a set of diameter at most $D$ is the union, over the at most $N(n,m)$ rank types, of the classes of each type in it, each of diameter at most $D$ by Lemma~\ref{lem:fibres}.
Lemma~\ref{lem:gluing-finite} with $C=1$, $T=16(n+m)r$, $\phi(D)=D$ and $N=N(n,m)$ completes the induction.
\end{proof}

\noindent Propositions~\ref{prop:dm-onegen} and~\ref{prop:dm-twoone} are the cases $n=1$ and $n=2$, $m\le1$, which they prove with explicit control: for one generator $\dM$ and $\di$ are bi-Lipschitz equivalent.
This equivalence fails in general. Whenever two non-isomorphic modules of $\Mnm$ share a rank invariant, scaling their grades keeps $\dM=0$ and sends $\di$ to infinity; $k(-(0,1))\oplus k(-(1,0))$ and $k[U_{\{(0,1),(1,0)\}}]\oplus k(-(1,1))$ do this in $\mathcal M_2^{3,1}$, with Betti numbers $(2,0)$ and $(3,1)$. And for $n,m\ge2$ over a field with at least four elements, the modules of Proposition~\ref{prop:twogen} are at matching distance $0$ and interleaving distance $t$, for every $t>0$ after scaling, so the quotient map $(\Mnm,\di)\to(\overline{\mathcal M}{}_2^{\,n,m},\overline d_{\mathcal M})$ is $1$-Lipschitz and fails to be a coarse equivalence.
Both sides have asymptotic dimension $2(n+m)$.
For $n,m\ge2$ over an uncountable field there is not even a coarse embedding of $(\Mnm,\di)$ into the quotient: the quotient is separable, each rank type being a Lipschitz image of a subset of a configuration space by Lemma~\ref{lem:fibres}, and a class of $(\Mnm,\di)$ containing the scaled families of Proposition~\ref{prop:twogen} has uncountable subsets at pairwise distance $t$ for every $t$, which a coarse embedding would carry to uncountable uniformly separated subsets of a separable space.
The constant in the Assouad--Nagata bound depends on $N(n,m)$ at each of the $2(n+m)$ steps of the induction.

\section{Discussion}\label{sec:disc}

Under a bound on the presentation the interleaving distance is tame at large scale: the asymptotic dimension and the Assouad--Nagata dimension are both the position count, $d(n+m)$ in $d$ parameters and $2(n+m)$ in two (Theorem~\ref{thm:position}).
The space is not doubling, even at a fixed fibered barcode (Corollary~\ref{cor:nodoubling}), and for $m\ge1$ it lies outside the proper setting of coarse geometry (Proposition~\ref{prop:notproper}).
The continuous moduli, which by scaling occur at every scale, are invisible at large scale---in the multiparameter analogue of $\asdim\D n=2n$.

\begin{remark}[Outside the proper setting]\label{rem:proper}
The operator-algebraic side of coarse geometry is formulated for proper metric spaces, usually of bounded geometry.
This is the setting of Yu's theorem that finite asymptotic dimension implies the coarse Baum--Connes conjecture~\cite{Yu1998} and of the coarse groupoid of Skandalis, Tu and Yu~\cite{SkandalisTuYu2002}.
By Proposition~\ref{prop:notproper} the spaces of this paper with $m\ge1$ and the diagram spaces $\D n$ fall outside that setting, and they have finite asymptotic dimension and embed coarsely into Hilbert space.
\end{remark}

The matching distance is computable~\cite{KerberLesnickOudot2019}, and it is the distance used in practice.
Under the matching distance, which sees only slices, the position count holds as well (Theorem~\ref{thm:matching}): slope-one slices pin the grid, and over a grid only finitely many rank invariants occur.
The isomorphism type of a grid representation is not determined by slices.

\section*{Declaration of AI-assisted technologies in the manuscript preparation}

During the preparation of this work, the authors used Claude (Anthropic) to draft and revise parts of the text, to verify citations, and for organization and editing. The authors reviewed and verified all content and take full responsibility for the content of the article.

\bibliographystyle{plain}
\bibliography{main}

\end{document}